\documentclass[]{siamart171218}

\usepackage{amsfonts}
\usepackage{subfig}
\usepackage{algorithm}
\usepackage{algorithmicx}

\newcommand{\R}{\mathbb{R}}

\newtheorem{remark}{Remark}[section]
\newtheorem{problem}{Problem}[section]

\title{On an inverse source problem for the full radiative transfer equation with incomplete data\thanks{Submitted to the editors DATE.
\funding{This work was supported by US Army Research Laboratory and US Army Research Office grant W911NF-19-1-0044. In addition, the effort of Nguyen was supported by research funds no.~FRG 111172 provided by The University of North Carolina at Charlotte.}}}

\author {Alexey V. Smirnov\thanks{Department of Mathematics and Statistics, University of North Carolina
Charlotte, Charlotte, NC, 28223
  (\email{asmirno2@uncc.edu}, \email{mklibanv@uncc.edu} (corresponding author), \email{loc.nguyen@uncc.edu}).}
\and Michael V. Klibanov\footnotemark[2]
\and  Loc H. Nguyen\footnotemark[2]}

\headers{Inverse Problem for Radiative Transfer Equation}{Alexey V. Smirnov, Michael V. Klibanov and Loc H. Nguyen}

\begin{document}

\maketitle

\begin{abstract}
A new numerical method to solve an inverse source problem for the radiative
transfer equation involving the absorption and scattering terms, with
incomplete data, is proposed. No restrictive assumption on those absorption
and scattering coefficients is imposed. The original inverse source problem
is reduced to boundary value problem for a system of coupled partial
differential equations of the first order. The unknown source function is
not a part of this system. Next, we write this system in the fully discrete
form of finite differences. That discrete problem is solved via the
quasi-reversibility method. We prove the existence and uniqueness of the
regularized solution. Especially, we prove the convergence of regularized
solutions to the exact one as the noise level in the data tends to zero
via a new discrete Carleman estimate. Numerical simulations demonstrate good performance of this method even when the data is highly noisy. 
\end{abstract}



\section{Introduction}

\label{sec 1}

The stationary radiative transfer equation (RTE) is commonly used in optics,
tomography, astrophysics, atmospheric science and remote sensing to describe
the propagation of the radiation field in media with absorbing, emitting and
scattering radiation. A significant number of studies is dedicated to the
recovery of the parameters of the observed objects from the measured data;
i.e., to the solutions of the inverse source problems (ISOPs) \cite{A07,Klo}
and coefficient inverse problems (CIPs) \cite{Bak,Silva}. A number of
inverse problems may be formulated, depending on the object's parameters of
one's interest.

The first reconstruction formula for the problem of the attenuated
tomography was obtained by Novikov \cite{Nov}. We also refer to \cite%
{Bal,Is,Nat2} for reconstruction formulae and as well as to \cite{N,Nat2}
for numerical results for the attenuated tomography with complete data and
with the scattering phase function $K\equiv 0$. Uniqueness and stability
results for similar ISOPs with complete data were obtained in \cite{Bal,Stef}%
. It was assumed in \cite{Bal} that $\left\vert K\right\vert $ is
sufficiently small. The assumption of \cite{Stef} is that functions $\sigma $
and $K$ belong to certain dense sets of some function spaces. The scattering
phase function $K$ is involved in RTE as the kernel of a certain integral
operator, the attenuation coefficient is $\sigma =\mu _{a}+\mu _{s}$, where $%
\mu _{a}$ and $\mu _{s}$ are the absorption and scattering coefficients
respectively, see Section 2.

In this paper, we propose a new numerical approach for the ISOP with limited
angle data for the stationary RTE and prove its convergence. This is \textit{%
the first publication,} in which a rigorously derived numerical method for
the ISOP for the RTE does not use any restrictive assumptions neither on $%
\mu _{a},$ nor on $\mu _{s},$ nor on $K$, except the smoothness and the
requirement that functions $\mu _{a}$ and $\mu _{s}$ are compactly
supported. Also, for the first time, a \textit{discrete} Carleman estimate
is applied here for the convergence analysis of an inverse problem. We note
that discrete Carleman estimates are very rare, unlike the continuous ones.
In addition, we prove the Lipschitz stability and uniqueness for our
statement of the ISOP.

Our method is based on the solution of an overdetermined boundary value
problem for a linear system of coupled integro-differential equations, in
which the unknown source function is not present. The solution of this
problem directly yields the solution of the desired ISOP. A similar idea was
recently used in \cite{KN}. However, unlike the current paper, a quite
restrictive condition\ $\sigma \equiv K\equiv 0$ is imposed in \cite{KN}.
The ISOP for the RTE with limited angle data has many applications in
optical imaging and tomography, such as bioluminescence tomography \cite%
{Jiang:07} and X-ray computerized tomography \cite{Lo86,Nat}.

The idea of our numerical method has roots in the Bukhgeim-Klibanov method
(BK) \cite{BK81}. BK was originally proposed in 1981 only for proofs of
global uniqueness and stability results for CIPs for PDEs, rather than for
numerical methods. BK is based on Carleman estimates. Since the current
paper is not a survey of BK, we refer here only to a few publications about
BK \cite{BK,Bell,K92,KT04,K13} and references cited therein. Currently, the
idea of BK is extensively used for constructions of globally convergent
numerical methods for CIPs for PDEs, see, e.g. \cite%
{Baud,KTsiam15,KK18,KLZ18}.

The second important element of our numerical method is the new orthonormal
basis in the space $L^{2}\left( a,b\right) ,\left( a,b\right) \subset 
\mathbb{R}$, which was recently introduced in \cite{Kjiip17}. This basis has
proven to be effective for numerical studies \cite{KK18,KLZ18,KN}. We use a
truncated Fourier series with respect to this basis. We estimate an optimal
number of terms of this series numerically and assume that this
approximation still satisfies the RTE, i.e. we work with an \textit{%
approximate mathematical model,} also, see Remark 4 at the end of Section 5.

We solve the above mentioned overdetermined boundary value problem by the
quasi-reversibility method (QRM), which is known to be effective to solve
overdetermined boundary value problems. We consider a fully discrete form of
our system, which is similar to what we use in the numerical tests. Next, we
establish a new discrete Carleman estimate and use it to prove uniqueness
and existence of the regularized solution for the QRM in the fully discrete
form, in which partial derivatives with respect to spatial variables are
written via finite differences. This Carleman estimate is also used to
establish the convergence rate of regularized solutions. Finally, we conduct
numerical testing for several different regimes of absorption and scattering
to show the method's potential for solving problems in real-world tomography.

The QRM was originally introduced by Lattes and Lions in 1969 \cite{LL}. We
also refer to, e.g. \cite{Bour1,Bour2,Is,KN} for this method. The second
author has shown in the survey paper \cite{KAP} that as long as a proper
Carleman estimate for an ill-posed problem for a linear PDE is available,
the convergent QRM can be constructed for this problem.

For brevity, we consider in this paper only the 2D case. The considerations
in the 3D case are similar. We state both forward and inverse problems in
Section 2. In Section 3 we derive the above mentioned over-determined
boundary value problem for a system of coupled partial differential
equations of the first order. To solve this problem, we apply the QRM by
stating a Minimization Problem. In Section 4 we introduce the fully discrete
version of the quasi-reversibility method to solve that problem. Next, we
derive a new discrete Carleman estimate. This estimate is used in Section 5
to prove the existence and uniqueness of the minimizer of the QRM and also
to establish the convergence rate of the minimizers to the exact solution as
the level of noise in the measured data tends to zero. Section 6 is devoted
to numerical studies. Everywhere below we work only with real-valued
functions.

\section{Statements of Forward and Inverse Problems}

\label{sec 2}

\begin{figure}
\begin{center}
\subfloat[The source/detector configuration of the problem in the case when the
source located at $\mathbf{x}_{\alpha }$ with $\left\vert \mathbf{x}_{\alpha
}\right\vert <R.$]{\includegraphics[width =.45\textwidth]{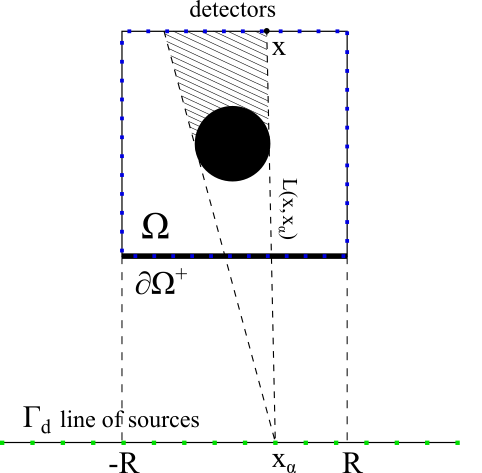}} \quad 
\subfloat[The source/detector configuration of the problem in the case when the
source located at $\mathbf{x}_{\alpha }$ with $\left\vert \mathbf{x}_{\alpha
}\right\vert >R.$]{\includegraphics[width =.45\textwidth]{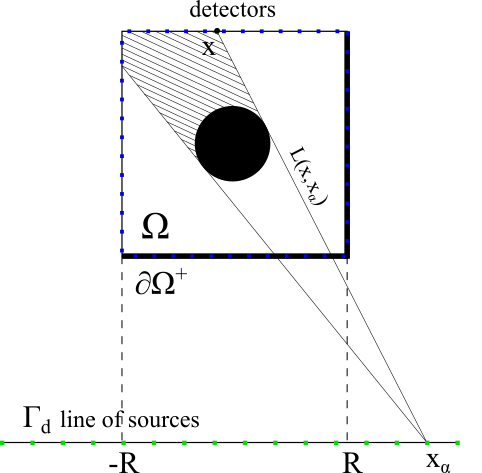}}
\end{center}
\caption{\textit{A schematic diagram of measurements for the 2D case. $L(\mathbf{x},\mathbf{x}_{\protect\alpha })$ is a straight line, connecting the detector $\mathbf{x}$ with the source $\mathbf{x}_{\protect\alpha }$.}}
\label{fig0}
\end{figure}

Let $\mathbf{x}=(x,y)$ denote an arbitrary point in $\mathbb{R}^{2}$. Let $%
a,b,d$ and $R$ be the positive numbers, where $1<a<b$ and $d\geq R$. Define
the rectangular domain $\Omega \subset \mathbb{R}^{2}$ (Figure \ref{fig0})
as 
\begin{equation}
\Omega =\{(x,y):-R<x<R,a<y<b\}.  \label{2.1}
\end{equation}

Let $\Gamma _{d}$ be the line with external sources 
\begin{equation}
\Gamma _{d}=\{\mathbf{x}_{\alpha }=(\alpha ,0):\alpha \in \lbrack -d,d]\}.
\label{2.2}
\end{equation}%
Let $u(\mathbf{x},\alpha )$ denotes the steady-state radiance at the point $%
\mathbf{x}$ generated by the external source located at $\mathbf{x}_{\alpha
}=(\alpha ,0)\in \Gamma _{d}$. Then, the function $u(\mathbf{x},\alpha )$
satisfies the following radiative transfer equation, see, e.g. \cite{Heino}%
\begin{multline}
    \nu (\mathbf{x},\alpha )\cdot \nabla _{\mathbf{x}}u(\mathbf{x},\alpha )+(\mu_{a}(\mathbf{x})+ \mu _{s}(\mathbf{x}))u(\mathbf{x},\alpha ) 
\\    =  
 \mu _{s}(\mathbf{%
x})\int_{\Gamma _{d}}K(\mathbf{x},\alpha ,\beta )u(\mathbf{x},\beta
)d\beta +f(\mathbf{x}) \quad \mbox{for all } \mathbf{x} \in \Omega.
\label{2.3}
\end{multline}
In the equation above, the function $f(\mathbf{x})\in L^{2}(\mathbb{R}^{2})$ is called the source
function while the functions $\mu _{a}(\mathbf{x})$, $\mu _{s}(\mathbf{x})\in
C^{1}\left( \mathbb{R}^{2}\right) $ denote the absorption and scattering
coefficients respectively. We assume that 
\begin{equation}
\mu _{a}(\mathbf{x})=\mu _{s}(\mathbf{x})=f\left( \mathbf{x}\right) =0 \quad \mbox{for all }
\mathbf{x}\in \mathbb{\R}^2 \setminus \Omega .  \label{2.4}
\end{equation}%
The function $K(\mathbf{x},\alpha ,\beta )\in C^{1}(\mathbb{R}^{2}\times
\lbrack -d,d]^{2})$ represents the so-called \textquotedblleft scattering
phase function". Scattering phase function is the probability density of a
particle scattering from $\nu (\mathbf{x},\beta )$-direction into $\nu (%
\mathbf{x},\alpha )$-direction. As the probability density, $K(\mathbf{x}%
,\alpha ,\beta )$ possesses the following properties, discussed in detail in 
\cite{Heino} 
\begin{equation}
K(\mathbf{x},\alpha ,\beta )\geq 0,\quad \int_{\Gamma _{d}}
\int_{\Gamma_{d}}K(\mathbf{x},\alpha ,\beta )d\alpha d\beta =1.  \label{2.5}
\end{equation}%
Finally, $\nu (\mathbf{x},\alpha )$ is the $\mathbb{R}^{2}-$vector, showing
the direction of particles propagating from the external source located at $%
\mathbf{x}_{\alpha }=(\alpha ,0)$ to $\mathbf{x,}$ 
\begin{equation}
\nu (\mathbf{x},\alpha )=\left( \frac{x-\alpha }{|\mathbf{x}-\mathbf{x}%
_{\alpha }|},\frac{y}{|\mathbf{x}-\mathbf{x}_{\alpha }|}\right) ,\quad x\in
\lbrack -R,R],\quad y\in \lbrack a,b],\quad \alpha \in \lbrack -d,d].
\label{2.6}
\end{equation}

For a fixed $\alpha $, let 
\begin{equation*}
\partial \Omega ^{+}=\left\{ \mathbf{x}\in \partial \Omega :\nu (\mathbf{x}%
,\alpha )\cdot \mathbf{n}(\mathbf{x})\leq 0\right\} ,
\end{equation*}%
where $\mathbf{n}(\mathbf{x})$ is the unit outward normal vector at $%
\partial \Omega $ at point $\mathbf{x}$. Assuming that all functions in
equation (\ref{2.3}), except ${u(\mathbf{x},\alpha)},$ are known
in $\Omega $, we formulate the following forward problem.
\begin{problem}[Forward Problem]
For each $\alpha \in \left[ -d,d\right],$ 
 find the function $\normalfont{u(\mathbf{x},\alpha)},$ satisfying
equation (\ref{2.3}) in the domain $\Omega $ as well as the
following boundary condition 
\begin{equation}
u(\mathbf{x},\alpha )=0 \quad \mbox{for all }\mathbf{x}\in \partial \Omega^{+}.%
\label{2.7}
\end{equation}%
\label{pro forward}
\end{problem}
In Appendix we prove existence and uniqueness of the solution of the
boundary value problem (\ref{2.3}), (\ref{2.7}) and; moreover, discuss a numerical
method to solve it. 
Conversely, assume now that the function $f(\mathbf{x})$ is
unknown and the information of $u(\mathbf{x}, \alpha)$ on $\partial \Omega$ is known. 
The main goal of this paper is to numerically solve the following inverse source problem:

\begin{problem}
[Inverse Source Problem] 
 Assume that equation \eqref{2.3} and conditions \eqref{2.4}, \eqref{2.5} hold. Also, let the vector $\nu(\mathbf{x}, \alpha)$ in \eqref{2.3} has the form \eqref{2.6}.
Reconstruct the function $f(\mathbf{x}),$ $\mathbf{x}%
\in \Omega,$ given the following boundary data 
\begin{equation}
F(\mathbf{x},\alpha )=u(\mathbf{x},\alpha ),\quad \mbox{for all } \mathbf{x}\in \partial
\Omega, \alpha \in [-d,d],  \label{2.8}
\end{equation}
where $\normalfont{u(\mathbf{x}, \alpha)}$ is the solution of Problem \ref{pro forward} and 
\begin{equation}
F(\mathbf{x},\alpha )=0\text{ }\textit{for}\text{ }\mathbf{x}\in \partial
\Omega ^{+}.  
\label{2.9}
\end{equation}
\label{ip}
\end{problem}

\begin{remark}
In the particular case when $\mu _{a}(\mathbf{x})\equiv \mu _{s}(\mathbf{x},\alpha )\equiv
0,$  this Inverse Source Problem is exactly the problem of X-ray
tomography with incomplete data, which was considered in \cite{KN}. However,
the main focus of this paper is to develop a numerical method for this problem allowing the presence of $\mu _{a},\mu
_{s}, K$. Especially, no technical condition is imposed on these interesting terms.
\end{remark}

\section{Numerical Method for the Inverse Source Problem}

\label{sec 3}

\subsection{An orthonormal basis in $L^2(-d, d)$}

First, we recall a special orthonormal basis in the space $L^{2}(-d,d)$,
which was introduced in \cite{Kjiip17}. For $\alpha \in \lbrack -d,d]$
consider the set of linearly independent functions $\left\{ \alpha
^{n-1}e^{\alpha }\right\} _{n=1}^{\infty }$. These functions form a complete
set in $L^{2}(-d,d)$. Applying the classical Gram-Schmidt orthonormalization
procedure to this set, we obtain the orthonormal basis $\{\Psi _{n}(\alpha
)\}_{n=1}^{\infty }$ in $L^{2}(-d,d)$. This basis has the following
properties \cite{Kjiip17}:
\begin{enumerate}
\item The functions $\Psi _{n}\in C^{1}[-d,d]$ and $\Psi _{n}^{\prime¶
}\left( \alpha \right) $ is not identically $0$, $\text{for all } \\ n=1,2,\dots $
\item 
$a_{nn}=1\text{ and }a_{mn}=0\text{ for all } m, n =1, 2, \dots$ such that $n<m,$ where
\begin{equation*}
a_{mn}=\int_{-d}^{d}\Psi _{n}^{\prime }(\alpha ),\Psi _{m}(\alpha
)d\alpha =
=\left\{ 
\begin{array}{ll}
1 & \mbox{if } m = m=n, \\ 
0 & \mbox{if } m > n.
m>n.%
\end{array}%
\right. 
\end{equation*}
\end{enumerate}
Item 2 implies that the matrix 
\begin{equation}
    M_{N}=(a_{mn})_{m,n=1}^{N}
    \label{matrix M}
\end{equation} \text{ is
invertible for all} $N=1,2,\dots.$

Hence, the function $u(\mathbf{x},\alpha )$ can be written as
the following Fourier series converging in $L^{2}(-d,d)$ for every point  $\mathbf{x} \in \overline \Omega$
\begin{equation*}
    u(\mathbf{x},\alpha ) = \sum_{n=1}^{\infty }u_{n}(\mathbf{x})\Psi_{n}(\alpha ) \quad \mbox{for all } \alpha \in [-d, d] 
    \end{equation*}
where
\[
    u_{n}(\mathbf{x}) =\int_{-d}^{d}u(\mathbf{x},\alpha )\Psi_{n}(\alpha) d\alpha.
\]
We approximate the function $u(\mathbf{x},\alpha )$ via the truncated Fourier series, and the same for $u_{\alpha }(\mathbf{x},\alpha ),$ 
\begin{align}
u(\mathbf{x},\alpha ) & \approx \sum_{n=1}^{N}u_{n}(\mathbf{x})\Psi _{n}(\alpha
),\quad \mathbf{x}\in \overline{\Omega },\quad \alpha \in \lbrack -d,d],
\label{3.1}
\\
u_{\alpha }(\mathbf{x},\alpha ) & \approx \sum_{n=1}^{N}u_{n}(\mathbf{x})\Psi
_{n}^{\prime }(\alpha ),\quad \mathbf{x}\in \overline{\Omega },\quad \alpha
\in \lbrack -d,d].  \label{3.2}
\end{align}%
where $N\geq 1$ is a certain integer, which is chosen numerically. We assume
that the truncated series (\ref{3.1}) satisfies equation (\ref{2.3}). In
addition, we assume that both sides of the equation resulting after the
substitution of (\ref{3.1}) in (\ref{2.3}) can be differentiated with
respect to the parameter $\alpha $ as in (\ref{3.2}). These assumptions form
our approximate mathematical model mentioned in Section \ref{sec 1}.

\subsection{A coupled system of first-order differential equations}

Just as in the first step of the above mentioned BK method \cite{BK81}, we
eliminate the unknown source function $f(\mathbf{x})$ from equation (\ref%
{2.3}) via the differentiation of that equation with respect to the
parameter $\alpha $ from which $f(\mathbf{x})$ does not depend. We obtain
\begin{multline}
\nu (x,y,\alpha )\cdot \nabla u_{\alpha }-\frac{y^{2}}{|\mathbf{x}-\mathbf{x%
}_{\alpha }|^{3}}u_{x}+\frac{(x-\alpha )y}{|\mathbf{x}-\mathbf{x}_{\alpha
}|^{3}} u_{y}+(\mu _{a}+\mu _{s})\left( \mathbf{x}\right) u_{\alpha } 
\\
-\mu _{s}(\mathbf{x})\int_{\Gamma _{d}}K_{\alpha }(\mathbf{x}%
,\alpha ,\beta )u(\mathbf{x},\beta )d\beta =0  \label{3.3}
\end{multline}
for all $\mathbf{x} = (x, y) \in \Omega.$
Multiplying equation (\ref{3.3}) by $|\mathbf{x}-\mathbf{x}_{\alpha }|/y$,
we obtain

\begin{align}
\begin{split}
& u_{y,\alpha }+\frac{x-\alpha }{y} u_{x,\alpha }+\frac{y}{|\mathbf{x}-%
\mathbf{x}_{\alpha }|^{2}}u_{x}+\frac{(x-\alpha )}{|\mathbf{x}-\mathbf{x}%
_{\alpha }|^{2}}u_{y} + \\
&+\frac{|\mathbf{x}-\mathbf{x}_{\alpha }|}{y}\left[ (\mu _{a}+\mu
_{s})\left( \mathbf{x}\right) u_{\alpha }-\mu _{s}(\mathbf{x}%
)\int_{\Gamma _{d}}K_{\alpha }(\mathbf{x},\alpha ,\beta )u(\mathbf{x}%
,\beta )d\beta \right] =0.  \label{3.4}
\end{split}%
\end{align}

Substituting representations (\ref{3.1}) and (\ref{3.2}) into equation (\ref%
{3.4}), multiplying the resulting equation by functions $\Psi _{m}(\alpha
),$ for each $m \in \{1,2,\dots ,N\}$  we obtain
\begin{multline}
\sum_{n=1}^{N}\frac{\partial u_{n}}{\partial y}\Psi _{n}^{\prime }(\alpha
)\Psi _{m}(\alpha )+\frac{x-\alpha }{y}\sum_{n=1}^{N}\frac{\partial u_{n}}{%
\partial x}\Psi _{n}^{\prime }(\alpha )\Psi _{m}(\alpha )
 \\
+\frac{y}{|\mathbf{x}-\mathbf{x}_{\alpha }|^{2}}\sum_{n=1}^{N}\frac{%
\partial u_{n}}{\partial x}\Psi _{n}(\alpha )\Psi _{m}(\alpha )+\frac{%
(x-\alpha )}{|\mathbf{x}-\mathbf{x}_{\alpha }|^{2}}\sum_{n=1}^{N}\frac{%
\partial u_{n}}{\partial y}\Psi _{n}(\alpha )\Psi _{m}(\alpha )
 \\
+\frac{|\mathbf{x}-\mathbf{x}_{\alpha }|}{y}\Psi _{m}(\alpha
)\sum_{n=1}^{N}\left[ (\mu _{a}+\mu _{s})\left( \mathbf{x}\right)
u_{n}\Psi _{n}^{\prime }(\alpha )\right] \\
-\frac{|\mathbf{x}-\mathbf{x}_{\alpha }|}{y}\Psi _{m}(\alpha
)\sum_{n=1}^{N}\left[\mu _{s}(\mathbf{x})\int_{\Gamma
_{d}}K_{\alpha }(\mathbf{x},\alpha ,\beta )u_{n}(\mathbf{x})\Psi _{n}\left(
\beta \right) d\beta \right] =0. 
 \label{3.5}
\end{multline}
Integrate equation (\ref{3.5}) with respect to $\alpha \in (-d,d).$
Recalling the definition of the matrix $M_{N}$ in \eqref{matrix M}, we obtain 
\begin{equation}
M_{N}U_{y}=\mathbf{A}U_{y}+\mathbf{B}U_{x}+\mathbf{C}U,\quad U(\mathbf{x}%
)=(u_{1}(\mathbf{x}),\dots ,u_{N}(\mathbf{x}))^{T},  \label{3.6}
\end{equation}%
Here \textbf{A},\textbf{B} and \textbf{C} are $N\times N$ matrices with the
following entries:
\begin{align}
(\mathbf{A})_{mn}& =\int_{\Gamma _{d}}\frac{(x-\alpha )}{|\mathbf{x}-%
\mathbf{x}_{\alpha }|^{2}}\Psi _{n}(\alpha )\Psi _{m}(\alpha )d\alpha ,
\label{3.7} \\
(\mathbf{B})_{mn}& =\int_{\Gamma _{d}}\left[ \frac{x-\alpha }{y}\Psi
_{n}^{\prime }(\alpha )\Psi _{m}(\alpha )+\frac{y}{|\mathbf{x}-\mathbf{x}%
_{\alpha }|^{2}}\Psi _{n}(\alpha )\Psi _{m}(\alpha )\right] d\alpha ,
\label{3.8} \\
\begin{split}
(\mathbf{C})_{mn}& =\int_{\Gamma _{d}}\frac{|\mathbf{x}-\mathbf{x}_{\alpha }|%
}{y} (\mu _{a}+\mu _{s})\left( \mathbf{x}\right) \Psi _{n}^{\prime }(\alpha
)\Psi _{m}(\alpha ) d\alpha \\
&-\int_{\Gamma _{d}}\frac{|\mathbf{x}-\mathbf{x}_{\alpha }|}{y} \mu _{s}(%
\mathbf{x})\left( \int_{\Gamma _{d}}K_{\alpha }(\mathbf{x},\alpha
,\beta )\Psi _{n}\left( \beta \right) d\beta\right) \Psi _{m}(\alpha )
d\alpha .  \label{3.9}
\end{split}%
\end{align}%
Everywhere below the norm of a matrix is the square root of the sum of
square norms of its entries. Since in the definition of the domain $\Omega $
the number $a>1$, the following estimates follow from (\ref{3.7})-(\ref%
{3.9}): 
\begin{equation*}
\max_{\mathbf{x}\in \overline{\Omega }}\big\vert\big\vert\mathbf{A}(%
\mathbf{x})\big\vert\big\vert\leq \frac{C_{0}}{a^{2}},\quad \max\limits_{%
\mathbf{x}\in \overline{\Omega }}\big\vert\big\vert\mathbf{B}(\mathbf{x})%
\big\vert\big\vert\leq \frac{C_{0}}{a},\quad \max\limits_{\mathbf{x}\in 
\overline{\Omega }}\big\vert\big\vert\mathbf{C}(\mathbf{x})\big\vert\big\vert%
\leq \frac{C_{0}}{a},
\end{equation*}%
where the number $C_{0}=C_{0}\left( R,d\right) >0$ depends only on the listed
parameters. Hence, there exists a sufficiently large number $%
a_{0}=a_{0}(N,R,d)>1$ such that for any $a>a_{0}$ the matrix $\tilde{\mathbf{%
A}}=M_{N}($Id$-M_{N}^{-1}\mathbf{A})$ is invertible. Everywhere below we
assume without further mentioning that $a>a_{0}.$ \newline

Denote $\mathbf{A}_{1} :=\tilde{\mathbf{A}}^{-1}\mathbf{B}$, $
\mathbf{A}_{2} := \tilde{\mathbf{A}}^{-1}\mathbf{C}$. Therefore,
equation (\ref{3.6}) is equivalent to 
\begin{equation}
U_{y}-\mathbf{A}_{1}U_{x}-\mathbf{A}_{2}U=0;\hspace{0.5em}\mathbf{A}_{1}=%
\mathbf{A}_{1}(x,y),\hspace{0.5em}\mathbf{A}_{2}=\mathbf{A}_{2}(x,y),\hspace{%
0.5em}(x,y)\in \Omega .  \label{3.10}
\end{equation}%
Using (\ref{2.8}) and (\ref{2.9}) we complement equation (\ref{3.10}) with
the following Dirichlet boundary condition 
\begin{equation}
U=F(x,y),\hspace{0.5em}\text{for}\hspace{0.5em}(x,y)\in \partial \Omega .
\label{3.11}
\end{equation}%
Thus, we have obtained a system of coupled linear differential
equations (\ref{3.10}) with the boundary condition (\ref{3.11}). The
solution $U(\mathbf{x})=(u_{1}(\mathbf{x}),\dots ,u_{N}(\mathbf{x}))^{T}$ of
the boundary value problem (\ref{3.10})--(\ref{3.11}) directly yields the desired
numerical solution to Problem \ref{ip} via the substitution of (\ref%
{3.1}) in (\ref{2.3}).

\subsection{The QRM for problem (\protect\ref{3.10})-- (\protect\ref{3.11})}

The problem (\ref{3.10})--(\ref{3.11}) is an overdetermined one. Indeed,
although equations (\ref{3.10}) are of the first order, the boundary
condition (\ref{3.11}) is given on the entire boundary $\partial \Omega .$
To find an approximate solution to this problem, we use the QRM, which, in
general works properly for overdetermined problems for PDEs. Thus, we
consider the following minimization problem for the Tikhonov-like functional 
$J_{\epsilon }$ with the regularization parameter $\epsilon \in \left(
0,1\right) :$ 
\begin{equation}
J_{\epsilon }(U)=\int_{\Omega }|U_{y}-\mathbf{A}_{1}U_{x}-\mathbf{A}%
_{2}U|^{2}dxdy+\epsilon \Vert \normalfont{U}\Vert _{H^{1}(\Omega )}^{2}.
\label{3.12}
\end{equation}%
When we say below that a vector function belongs to a Hilbert space, we mean
that each of its components belongs to this space and its norm is the square root of the sum of norms in that space of its components. \newline

\begin{problem}[Minimization Problem]
 Minimize the functional $J_{\epsilon
}$ on the set of $\mathit{N}$-dimensional vector valued functions $U\in
H^{1}(\Omega )$ satisfying boundary condition (\ref{3.11}). \label{minimize}
\label{pro min}
\end{problem}

\section{The Fully Discrete Form of the QRM}

\label{sec 4}

To solve Problem \ref{pro min}, we write $U_{x},U_{y}$ in the functional $%
J_{\epsilon }(U)$ in its finite difference version and minimize it with respect to
values of the vector function $U$ at grid points. Hence, we formulate the
QRM in this section in the fully discrete form of finite differences. We
prove existence and uniqueness of the minimizer and establish convergence
rate of minimizers to the exact solution, which is also written via finite
differences.

\subsection{The fully discrete form of functional (\protect\ref{3.12})}

Consider the following uniform 2-dimensional grid points on $\overline{\Omega }$ whose $x$ and $y$ coordinates are given by
\begin{align}
-R& =x_{0}<x_{1}<\dots <x_{M_x}=R,\hspace{0.5em} x_{i+1}-x_{i}=h_{x},\hspace{%
0.5em}\forall i\in \{0,1,\dots ,M_x-1\},  \label{4.1} \\
a& =y_{0}<y_{1}<\dots <y_{M_{y}}=b,\hspace{0.5em} y_{j+1}-y_{j}=h_{y},%
\hspace{0.5em}\forall j\in \{0,1,\dots ,M_{y}-1\}.  \label{4.2}
\end{align}%
Denote $h=(h_{x},h_{y}).$ We define the discrete set $\Omega ^{h}\hspace{%
0.5em}$as 
\begin{align*}
\Omega ^{h}& =\{(x,y):\left\{ \left( x_{i},y_{j}\right) \right\} ,i\in
\{1,\dots ,M_{x}-1\},\hspace{0.5em}j\in \{1,\dots ,M_{y}-1\}\}, \\
\partial \Omega ^{h}& =\{(x,y):\left\{ \left( x_{i},y_{j}\right) \right\} 
\hspace{0.5em}\text{for}\hspace{0.5em}i=0,M_{x},\hspace{0.5em}j=0,M_{y}\}, \\
\overline{\Omega }^{h}& =\Omega ^{h}\cup \partial \Omega ^{h}.
\end{align*}%
For any $N-$D matrix $Q(x,y)\in C(\overline{\Omega })$ we introduce the
following notations 
\begin{align}
\begin{split}
\mathbf{Q}_{i,j}^{h}& =Q(x_{i},y_{j}),\hspace{0.5em}i\in \{1,\dots
,M_{x}-1\},\hspace{0.5em}j\in \{1,\dots ,M_{y}-1\}, \\
\widetilde{\mathbf{Q}}_{i,j}^{h}& =Q(x_{i},y_{j}),\hspace{0.5em}i\in
\{0,\dots ,M_{x}\},\hspace{0.5em}j\in \{0,\dots ,M_{y}\}, \\
\mathbf{Q}^{h} \hspace{0.375em}& =\left\{ \mathbf{Q}_{i,j}^{h}\right\}
_{i,j=1}^{M_{x}-1,M_{y}-1}\hspace{0.5em}\text{is an }(M_{x}-1)\times
(M_{y}-1)\text{ matrix}, \\
\widetilde{\mathbf{Q}}^{h} \hspace{0.375em}& =\left\{ \mathbf{Q}%
_{i,j}^{h}\right\} _{i,j=0}^{M_{x},M_{y}}\hspace{0.5em}\text{is an }%
(M_{x}+1)\times (M_{y}+1)\text{ matrix}.  \label{4.3}
\end{split}%
\end{align}

Note that the matrix $\mathbf{Q}^{h},$ in contrast to $\widetilde{\mathbf{Q}}%
^{h},$ does not include boundary terms of the form 
\begin{equation*}
\mathbf{Q}_{0,j}^{h}=Q(-R,y_{j}),\hspace{0.5em}\mathbf{Q}%
_{i,0}^{h}=Q(x_{i},a),\hspace{0.5em}\mathbf{Q}_{M_{x},j}^{h}=Q(R,y_{j}),%
\hspace{0.5em}\mathbf{Q}_{i,M_{y}}^{h}=Q(x_{i},b).
\end{equation*}

Recall the forward finite difference formulae for the vector function $%
\mathbf{Q}^{h}$: 
\begin{align}
& (\mathbf{U}^{h})_{x}^{\prime }=\left\{ (\mathbf{U}_{i,j}^{h})_{x}^{\prime
}\right\} _{i,j=0}^{M_{x}-1,M_{y}-1},(\mathbf{U}_{i,j}^{h})_{x}^{\prime }=%
\frac{\mathbf{U}_{i+1,j}^{h}-\mathbf{U}_{i,j}^{h}}{h_{x}},\hspace{0.5em}i\in
\{0,\dots ,M_{x}-1\},  \label{4.4} \\
& (\mathbf{U}^{h})_{y}^{\prime }=\left\{(\mathbf{U}_{i,j}^{h})_{y}^{\prime }
\right\}_{i,j=0}^{M_{x}-1,M_{y}-1}, (\mathbf{U}_{i,j}^{h})_{y}^{\prime }=%
\frac{\mathbf{U}_{i,j+1}^{h}-\mathbf{U}_{i,j}^{h}}{h_{y}},\hspace{0.5em}j\in
\{0,\dots ,M_{y}-1\}.  \label{4.5}
\end{align}
Hence, we obtain the following finite difference analog of (\ref{3.10})--(\ref{3.11}) 
\begin{align}
L^{h}\left( \mathbf{U}^{h}\right) & =(\mathbf{U}^{h})_{y}^{\prime }-\mathbf{A%
}_{1}^{h}(\mathbf{U}^{h})_{x}^{\prime }+\mathbf{A}_{2}^{h}\mathbf{U}^{h}=0;%
\hspace{0.5em}\mathbf{A}_{1}^{h}=\mathbf{A}_{1},\hspace{0.5em}\mathbf{A}%
_{2}^{h}=\mathbf{A}_{2}\hspace{0.5em}\text{in}\hspace{0.5em}\Omega ^{h},
\label{4.6} \\
\widetilde{\mathbf{U}}^{h}& =\mathbf{F}^{h}\hspace{0.5em}\text{on}\hspace{%
0.5em}\partial \Omega ^{h}.  \label{4.7}
\end{align}%
where the boundary matrix $\mathbf{F}^{h}$ is defined using the values of
the matrix $F(x,y)$ on the grid (\ref{4.1}), (\ref{4.2}). We define the
following discrete functional spaces for matrices $\mathbf{Q}^{h}$, $%
\widetilde{\mathbf{Q}}^{h}$:
\begin{equation*}
L^{2,h}(\Omega ^{h}) =\Big\{ \mathbf{Q}^{h}:\Vert \mathbf{Q}^{h}\Vert
_{L^{2,h}(\Omega
^{h})}^{2}=h_{y}h_{x}\sum_{j=1}^{M_{y}-1}\sum_{i=1}^{M_{x}-1}[\mathbf{Q}%
_{i,j}^{h}]^{2}<\infty \Big\}, 
\end{equation*}
and
\begin{multline*}
H^{1,h}(\Omega ^{h}) =\Big\{ 
\mathbf{Q}^{h}:\Vert \mathbf{Q}^{h}\Vert _{H^{1,h}(\Omega ^{h})}^{2} 
\\ 
=h_{y}h_{x}\sum_{i=1}^{M_{x}-1}\sum_{j=1}^{M_{y}-1}([(%
\mathbf{Q}_{i,j}^{h})_{x}^{\prime }]^{2}+[(\mathbf{Q}_{i,j}^{h})_{y}^{\prime
}]^{2}+[\mathbf{Q}_{i,j}^{h}]^{2})<\infty%
\Big\}.
\end{multline*}
We define the inner products in these spaces in the obvious manner and
denote them as $\big(\cdot ,\cdot \big)$ and $\big[\cdot ,\cdot \big]$ for $%
L^{2,h}(\Omega ^{h})$ and $H^{1,h}(\Omega ^{h})$ respectively.

\begin{remark}
  Here and everywhere below if a matrix $%
{\mathbf{Q}^h}$ is defined as in {(\ref{4.3})}, ~then~$%
{\widetilde{\mathbf{Q}}^h}$ denotes the matrix $%
{\mathbf{Q}^h}$, complemented by boundary conditions at $%
{\partial\Omega^h}$. 
\end{remark}

\begin{remark}
 Below we fix the number $%
h_{1}\in (0,1)$ and restrict $h_{x}$ from the below as $h_{x}\in \left[
h_{1},1\right) .$ However, we do not restrict from the below $h_{y}>0$ by a
positive constant. Then it follows from (\ref{4.5}) that there
exists a constant $B_{h_{1}}>0$ depending only on $h_{1}$ such that if $\mathbf{Q}_{0,j}^{h}=\mathbf{Q}_{M_{x},j}^{h}=0;j=1,...,M_{y}-1,$ then 
\begin{equation}
 \Vert (\mathbf{Q}^h)^{\prime}_x \Vert_{L^{2,h}(\Omega^h)}^2 \leq\quad
B_{h_1} \Vert \mathbf{Q}^h \Vert_{L^{2,h}(\Omega^h)}^2, \quad \forall h
\in [h_1,1], \quad \forall \mathbf{Q}^h :
\normalfont{\widetilde{\mathbf{Q}}^h} \in H^{1,h}(\Omega^h).  \label{4.8}
\end{equation}
\end{remark}

The fully discrete QRM applied to problem (\ref{4.6})-(\ref{4.7}) leads to
the following discrete version of the above Minimization Problem: 
\begin{problem}[Discrete Minimization Problem]
Minimize the functional 
\begin{equation}
J_{\epsilon }^{h}(\normalfont{\widetilde{\mathbf{U}}^h})=\Vert %
\normalfont{({\mathbf{U}}^h)^{\prime}_y}-\mathbf{A}_{1}^{h}\normalfont{({%
\mathbf{U}}^h)^{\prime}_x}-\mathbf{A}_{2}^{h}\normalfont{{\mathbf{U}}^h}\Vert
_{L^{2,h}({\Omega }^{h})}^{2}+\epsilon \Vert \normalfont{\mathbf{U}^h}\Vert
_{H^{1,h}(\Omega ^{h})}^{2}  \label{4.9}
\end{equation}%
on the set of matrices $\normalfont{\widetilde{\mathbf{U}}}^{h}$, satisfying
the boundary condition \eqref{4.7}. 
\end{problem}
The minimizer of $J_{\epsilon }^{h}(\normalfont{\widetilde{\mathbf{U}}^h})$
satisfying boundary condition (\ref{4.7}) is called the \textit{regularized
solution} of the problem (\ref{4.6})--(\ref{4.7}).

\subsection{A discrete Carleman estimate}

We now derive a discrete Carleman estimate for the finite difference version
of the differential operator $d/dy$. Consider a uniform partition of the
interval $(a,b)\subset \mathbb{R}$ of the real line into $M$ subintervals
with the grid step size $h_{y}$, 
\begin{equation}
a=y_{0}<y_{1}<\ldots <y_{M}=b,\quad y_{j+1}-y_{j}=h_{y}, \quad j\in \{0,1,\ldots
,M-1\}.  \label{4.10}
\end{equation}%
Following the book \cite{Sam}, for any discrete function $u$ defined on this
grid denote 
$u_{j}=u(y_{j})$ and define both its forward $u_{j}^{\prime }$ and backward $%
\overline{u}_{j}^{\prime }$ finite difference derivatives, which are the
finite difference analogs of the differential operator $d/dy,$ as 
\begin{equation}
u_{j}^{\prime }=\frac{(u_{j+1}-u_{j})}{h_{y}},j\in \{0,\ldots ,M-1\},\quad 
\overline{u}_{j}^{\prime }=\frac{(u_{j}-u_{j-1})}{h_{y}},\forall j\in
\{1,\ldots ,M\}.  \label{4.11}
\end{equation}

\begin{lemma}
\textit{For any discrete function $w$, defined on the grid} \normalfont{(\ref{4.10})}
\textit{the following inequality holds: }
\begin{equation*}
-2h_{y}\sum_{j=1}^{M-1}w_{j}w_{j}^{\prime }\geq -(w_{M}^{2}-w_{1}^{2}).
\label{lem1}
\end{equation*}
\end{lemma}

\textbf{Proof.} Using the summation by parts formula for the discrete
function $w$ \cite{Sam}, we obtain 
\begin{equation*}
h_{y}\sum_{j=1}^{M-1}w_{j}w_{j}^{\prime
}=(w_{M}^{2}-w_{0}w_{1})-h_{y}\sum_{i=1}^{M}w_{i}\overline{w_{i}}^{\prime }.
\end{equation*}%
Next, 
\begin{align*}
 h_{y}\sum_{i=1}^{M}w_{i}\overline{w_{i}}^{\prime
}
    &=h_{y}\sum_{j=0}^{M-1}w_{j+1}w_{j}^{\prime
}
    =h_{y}\sum_{j=0}^{M-1}(w_{j}+w_{j}^{\prime }h_{y})w_{j}^{\prime } 
\\
& =h_y (w_{0}+hw_{0}^{\prime })w_{0}^{\prime
}
+h_{y}\sum_{j=1}^{M-1}(w_{j}+w_{j}^{\prime }h_{y})w_{j}^{\prime
}
\\
&=w_{1}(w_{1}-w_{0})+h_{y}\sum_{j=1}^{M-1}(w_{j}+w_{j}^{\prime
}h_{y})w_{j}^{\prime }
\end{align*}%
Combining all equalities written above, we obtain 
\begin{align*}
h_{y}\sum_{j=1}^{M-1}w_{j}w_{j}^{\prime } 
&=  \left(
w_{M}^{2}-w_{0} w_1\right) - (w_1^2 -w_0 w_1) -h_{y}\sum_{j=1}^{M-1}\left( w_{j}+w_{j}^{\prime
}h_{y}\right) w_{j}^{\prime } \\
&= \left( w_{M}^{2}-w_{1}^{2}\right)
-h_{y}\sum_{j=1}^{M-1}w_{j}w_{j}^{\prime
}-h_{y}^{2}\sum_{j=1}^{M-1}(w_{j}^{\prime })^{2}.
\end{align*}%
Hence,%
\begin{equation*}
-2h_{y}\sum_{j=1}^{M-1}w_{j}w_{j}^{\prime
}=-(w_{M}^{2}-w_{1}^{2})+h_{y}^{2}\sum_{j=1}^{M-1}(w_{j}^{\prime })^{2}\geq-
\left( w_{M}^{2}-w_{1}^{2}\right) .\text{ }\square
\end{equation*}

\begin{theorem}[A discrete Carleman estimate]
For any positive number $\lambda >0$, the following discrete Carleman
estimate holds for any discrete function $u$, defined on the grid (\ref{4.10}) 
\begin{multline*}
h_{y}\sum_{j=1}^{M_{y}-1}e^{2\lambda y_{j}}\left( u_{j}^{\prime }\right)
^{2}\geq h_{y}\sum_{j=1}^{M_{y}-1}  \left( \frac{1- e^{-%
\lambda h_{y}}}{h_{y}}\right) ^{2}e^{2\lambda y_{j}}u_{j}^{2}
 \\
+ 2e^{-\lambda h_{y}}  \left( \frac{1- e^{-%
\lambda h_{y}}}{h_{y}}\right) (e^{2\lambda y_{1}}u_{1}^{2}- %
e^{2\lambda y_{M}}u_{M}^{2}).
\end{multline*}%
\label{thm1}
\end{theorem}

\textbf{Proof}. For each $j$, we define 
\begin{equation}
w_{j}=e^{\lambda y_{j}}u_{j}, \quad u_{j}=e^{-\lambda y_{j}}w_{j}.
\label{4.12}
\end{equation}%
Hence, according to (\ref{4.11}), the forward difference derivative of the
function $u$ at $y_{j}$ is

\begin{align*}
u_{j}^{\prime } &=\frac{e^{-\lambda (y_{j}+h_{y})}w_{j+1}-e^{-\lambda
y_{j}}w_{j}}{h_{y}}
 =\frac{e^{-\lambda y_{j}}\left( e^{-\lambda
h_{y}}w_{j+1}-w_{j}\right) }{h_{y}} \\
&=e^{-\lambda y_{j}}\left( e^{-\lambda h_{y}}\frac{w_{j+1}-w_{j}}{h_{y}} + 
\frac{e^{-\lambda h_{y}}w_{j}-w_{j}}{h_{y}}\right) 
\\
&=e^{-\lambda y_{j}}\left(
w_{j}^{\prime }e^{-\lambda h_{y}}-\frac{1-e^{-\lambda h_{y}}}{h_{y}}%
w_{j}\right) .
\end{align*}
Hence, we have for each $j=1,...,M_{y}-1:$ 
\begin{align*}
e^{2\lambda y_{j}}(u_{j}^{\prime })^{2}&=\left( w_{j}^{\prime }e^{-\lambda
h_{y}}-\frac{1-e^{-\lambda h_{y}}}{h_{y}}w_{j}\right) ^{2} \\
&\geq \left( \frac{1-e^{-\lambda h_{y}}}{h_{y}}\right) ^{2}\left( w_{j}\right)
^{2}-2e^{-\lambda h_{y}}\frac{\left( 1-e^{-\lambda h_{y}}\right) }{h_{y}}%
w_{j}^{\prime }w_{j}.
\end{align*}
As a result,
\begin{equation*}
h_{y}\sum_{j=1}^{M_{y}-1}e^{2\lambda y_{j}}(u_{j}^{\prime })^{2}\geq
h_{y}\sum_{j=1}^{M_{y}-1}\left( \frac{1- e^{-\lambda
h_{y}}}{h_{y}}\right) ^{2}w_{j}^{2}-2e^{-\lambda h_{y}}%
\frac{(1- e^{-\lambda h_{y}})}{h_{y}}h_{y}\sum_{j=1}^{M%
-1}w_{j}^{\prime }w_{j}.
\end{equation*}%
Applying Lemma \ref{lem1} to the second term in the right hand side, we obtain 
\begin{equation*}
\hspace{1em}h_{y}\sum_{j=1}^{M_{y}-1}e^{2\lambda y_{j}}(u_{j}^{\prime
})^{2}\geq h_{y}\sum_{j=1}^{M_{y}-1}\left( \frac{1- e^{-%
\lambda h_{y}}}{h_{y}}\right) ^{2}w_{j}^{2}+\ 2e^{-\lambda
h_{y}}\left( \frac{1- e^{-\lambda h_{y}}}{h_{y}}\right)
(w_{1}^{2}- w_{M}^{2}).
\end{equation*}%
The statement of Theorem \ref{thm1} follows from this estimate and %
\eqref{4.12}. $\square $

\begin{lemma}
Let $u$ be a discrete function, defined on the grid \eqref{4.10}, such that $%
u_{M}=0$. 
Then for any two numbers $\lambda ,h_{y}>0$ such that $\lambda
h_{y}<1$ the following inequality holds 
\begin{equation}
h_{y}\sum_{j=1}^{M_{y}-1}e^{2\lambda y_{j}}\left( u_{j}^{\prime }\right)
^{2}\geq \frac{\lambda ^{2}}{4}h_{y}\sum_{j=1}^{M_{y}-1}e^{2\lambda
y_{j}}u_{j}^{2}.  \label{4.13}
\end{equation}%
\label{lem2}
\end{lemma}

\textbf{Proof.} By Taylor formula 
\begin{equation*}
e^{-\lambda h_{y}}=1-\lambda h_{y}+\frac{e^{-\xi }}{2}\left( \lambda
h_{y}\right) ^{2}=1-\lambda h_{y}\left( 1-\frac{e^{-\xi }}{2}\lambda
h_{y}\right) ,
\end{equation*}%
where $\xi \in \left( 0,\lambda h_{y}\right) $ is a certain number. Hence, $%
1-e^{-\lambda h_{y}}\geq \lambda h_{y}/2.$ Hence,  
\begin{equation*}
\left( \frac{1- e^{-\lambda h_{y}}}{h_{y}}\right)
^{2}\geq \frac{\lambda ^{2}}{4}.
\end{equation*}%
Therefore, using Theorem \ref{thm1}, we obtain (\ref{4.13}). $\square $ 

\begin{remark}
 This lemma is a discrete analog of the Carleman
estimate in \cite[Lemma 4.1]{KN} for the continuous case of the operator $d/dy.$ 
\end{remark}

\section{Convergence Analysis}

\label{sec 5}

\subsection{Existence of the solution of the Discrete Minimization Problem}

\begin{theorem}
Assume that there exists a matrix $\mathbf{G}^{h}\in H^{1,h}(\Omega
^{h})$ such that $\widetilde{\mathbf{G}}^{h}\mid _{\partial \Omega
^{h}}=\mathbf{F}^{h}.$ Then for each $\epsilon >0,$ there
exists unique minimizer $\mathbf{U}_{\min ,\epsilon }^{h}\in {H}%
^{1,h}(\Omega ^{h})$ of the functional (\ref{4.9}) satisfying boundary
condition (\ref{4.7}). \label{thm2}
\end{theorem}

\textbf{Proof.} Let $H_{0}^{1,h}(\Omega ^{h})$ be the subspace of the space $%
H^{1,h}(\Omega ^{h})$ consisting on such matrices $\widetilde{\mathbf{Q}}%
^{h}\in H^{1,h}(\Omega ^{h})$ that $\widetilde{\mathbf{Q}}^{h}\mid
_{\partial \Omega ^{h}}=0.$ Recalling notation (\ref{4.6}) for the operator $%
L^{h},$ we rewrite the functional $J_{\epsilon }^{h}(\normalfont{\widetilde{%
\mathbf{U}}^h})$ in the following form:%
\begin{align}
J_{\epsilon }^{h}(\mathbf{W}^{h}) &=\left\Vert L^{h}\left( \mathbf{W}%
^{h}\right) +L^{h}\left( \mathbf{G}^{h}\right) \right\Vert _{L^{2,h}\left(
\Omega ^{h}\right) }^{2}+\epsilon \left\Vert \mathbf{W}^{h}+\mathbf{G}%
^{h}\right\Vert _{H^{1,h}(\Omega ^{h})}^{2},  \label{5.1} \\
\widetilde{\mathbf{W}}^{h}\mid _{\partial \Omega } &= (\widetilde{\mathbf{U}}%
^{h}-\widetilde{\mathbf{G}}^{h})\mid _{\partial \Omega } = 0.  \label{5.2}
\end{align}%
Thus, in order to work with zero boundary condition in (\ref{5.1})-(\ref{5.2}%
), we consider the function $\mathbf{W}^{h}=\mathbf{U}^{h}-\mathbf{G}^{h}$
instead of $\mathbf{U}^{h}.$ \newline

Let $\mathbf{W}_{\min ,\epsilon }^{h}$ with $\widetilde{\mathbf{W}}_{\min
,\epsilon }^{h}\in H_{0}^{1,h}(\Omega ^{h})$ be any minimizer of functional (%
\ref{5.1}). By the variational principle the following identity holds for
all $\widetilde{\mathbf{P}}^{h}\in H_{0}^{1,h}(\Omega ^{h}):$ 
\begin{multline}
\left( L^{h}\left( \mathbf{W}_{\min ,\epsilon }^{h}\right) ,L^{h}\left( 
\mathbf{P}^{h}\right) \right) +\epsilon \left[ \mathbf{W}_{\min ,\epsilon
}^{h},\mathbf{P}^{h}\right] =-\left( L^{h}\left( \mathbf{G}^{h}\right)
,L^{h}\left( \mathbf{P}^{h}\right) \right) 
\\
 -\epsilon \left[ \mathbf{G}^{h},\mathbf{P}^{h}\right] .
\label{5.3}
\end{multline}%
The left hand side of the identity (\ref{5.3}) generates a new scalar
product $\left\{ \cdot ,\cdot \right\} $ in the subspace $H_{0}^{1,h}(\Omega
^{h}).$ Consider the corresponding norm $\left\{ \cdot \right\} ^{2},$%
\begin{equation}
\left\{ \mathbf{Q}^{h}\right\} ^{2}=\left\Vert L^{h}\left( \mathbf{Q}%
^{h}\right) \right\Vert _{L^{2,h}\left( \Omega ^{h}\right) }^{2}+\epsilon
\left\Vert \mathbf{Q}^{h}\right\Vert _{H^{1,h}(\Omega ^{h})}^{2},\forall 
\widetilde{\mathbf{Q}}^{h}\in H_{0}^{1,h}(\Omega ^{h}).  \label{5.4}
\end{equation}%
Obviously, there exists a certain constant $C_{1}=C_{1}\left( L^{h},h,\Omega ^{h},\epsilon
\right) >0,$ which depends only on listed parameters such that  
\begin{equation}
\epsilon \left\Vert \mathbf{Q}^{h}\right\Vert _{H^{1,h}(\Omega
^{h})}^{2}\leq \left\{ \mathbf{Q}^{h}\right\} ^{2}\leq C_{1}\left\Vert 
\mathbf{Q}^{h}\right\Vert _{H^{1,h}(\Omega ^{h})}^{2},\forall \widetilde{%
\mathbf{Q}}^{h}\in H_{0}^{1,h}(\Omega ^{h}),  \label{5.5}
\end{equation}%
Below $C_{1}$ denotes
different positive numbers depending on the same parameters.

Hence, norms $\left\{ \mathbf{Q}^{h}\right\} $ and $\left\Vert \mathbf{Q}%
^{h}\right\Vert _{H^{1,h}(\Omega ^{h})}$ are equivalent for $\widetilde{%
\mathbf{Q}}^{h}\in H_{0}^{1,h}(\Omega ^{h}).$ Therefore, (\ref{5.3}) is
equivalent with 
\begin{equation}
\left\{ \widetilde{\mathbf{W}}_{\min ,\mu }^{h},\mathbf{\mathbf{P}}%
^{h}\right\} =-\left( L^{h}\left( \mathbf{G}^{h}\right) ,L^{h}\left( \mathbf{%
P}^{h}\right) \right) -\epsilon \left[ \mathbf{G}^{h},\mathbf{P}^{h}\right]
,\forall \widetilde{\mathbf{P}}^{h}\in H_{0}^{1,h}(\Omega ^{h}).  \label{5.6}
\end{equation}%
Using the Cauchy-Schwarz inequality, (\ref{5.4}) and (\ref{5.5}), we obtain 
\begin{equation*}
\left\vert -\left( L^{h}\left( \mathbf{G}^{h}\right) ,L^{h}\left( \mathbf{P}%
^{h}\right) \right) -\epsilon \left[ \mathbf{G}^{h},\mathbf{P}^{h}\right]
\right\vert \leq C_{1}\left\{ \mathbf{G}^{h}\right\} \left\{ \mathbf{P}%
^{h}\right\}, \quad \mbox{for all } \widetilde{\mathbf{P}}^{h}\in H_{0}^{1,h}(\Omega ^{h}).
\end{equation*}%
Hence, the right hand side of (\ref{5.6}) can be considered as a bounded
linear functional mapping the space $H_{0}^{1,h}(\Omega ^{h})$ in $\mathbb{R}
$. Since the regular norm in $H_{0}^{1,h}(\Omega ^{h})$ is equivalent with
the norm generated by new scalar product $\left\{ \cdot ,\cdot \right\} ,$
then Riesz representation theorem implies that there exists unique matrix $%
\widetilde{\boldsymbol{\Phi }}^{h}\in H_{0}^{1,h}(\Omega ^{h})$ such that 
\begin{equation*}
\left\{ \widetilde{\mathbf{W}}_{\min ,\epsilon }^{h},\widetilde{\mathbf{P}}%
^{h}\right\} =\left\{ \widetilde{\boldsymbol{\Phi }}^{h},\widetilde{\mathbf{P%
}}^{h}\right\} ,\quad \mbox{for all }\widetilde{\mathbf{P}}^{h}\in
H_{0}^{1,h}\left( \Omega ^{h}\right) .
\end{equation*}%
Therefore, $\widetilde{\mathbf{W}}_{\min ,\epsilon }^{h}=\widetilde{\boldsymbol{%
\Phi }}^{h}\hspace{0.5em}$and$\hspace{0.5em}\mathbf{W}_{\min ,\epsilon }^{h}=%
\boldsymbol{\Phi }^{h}.$ Finally, the matrix $\mathbf{U}_{\min ,\epsilon
}^{h}=\mathbf{W}_{\min ,\epsilon }^{h}+\mathbf{G}^{h}$ is the unique
minimizer claimed by this theorem. $\square $

\subsection{Convergence rate of regularized solutions, Lipschitz stability
and uniqueness}

The minimizer $\widetilde{\mathbf{U}}_{\min ,\epsilon }^{h}$ is called the 
\textit{regularized solution} of problem (\ref{4.6}), (\ref{4.7}). In this
section, we establish the convergence rate of regularized solutions to the
exact one when the noise in the data tends to zero. \ In addition, we
establish Lipschitz stability estimate and uniqueness for the problem (\ref{4.6}%
), (\ref{4.7}). 

Let a matrix $\mathbf{P}^{h}\in L^{2,h}\left( \Omega ^{h}\right) .$ Denote 
\begin{equation}
\left\Vert \mathbf{P}^{h}e^{\lambda y}\right\Vert _{L^{2,h}\left( \Omega
^{h}\right)
}^{2}=h_{y}h_{x}\sum_{j=1}^{M_{y}-1}\sum_{i=1}^{M_{x}-1}%
\left( \mathbf{P}_{i,j}^{h}\right) ^{2}e^{2\lambda y_{j}}.  \label{5.7}
\end{equation}%
Hence, by Lemma \ref{lem2} for all $\lambda h_{y}\in \left( 0,1\right) $ and for all $%
\widetilde{\mathbf{P}}^{h}\in H_{0}^{1,h}\left( \Omega ^{h}\right) $ 
\begin{equation}
\left\Vert \left( \mathbf{P}^{h}\right) _{y}^{\prime }e^{\lambda
y}\right\Vert _{L^{2,h}\left( \Omega ^{h}\right) }^{2}\geq \frac{\lambda
^{2}}{8}\left\Vert \mathbf{P}^{h}e^{\lambda y}\right\Vert _{L^{2,h}\left( \Omega
^{h}\right) }^{2}.  \label{5.8}
\end{equation}

Let $\mathbf{U}^{\ast h}\in H^{1,h}\left( \Omega ^{h}\right) $ be the exact
solution of problem (\ref{4.6}), (\ref{4.7}) with the exact boundary data $%
\mathbf{F}^{\ast ,h}.$ We assume that there exists an exact matrix $\mathbf{G%
}^{\ast h}$ such that 
\begin{equation}
\mathbf{G}^{\ast h}\in H^{1,h}(\Omega ^{h}),\text{ }\widetilde{\mathbf{G}}%
^{\ast ,h}\mid _{\partial \Omega ^{h}}=\mathbf{F}^{\ast ,h}.  \label{5.9}
\end{equation}%
As to the boundary data $\mathbf{F}^{h},$ we assume, as in Theorem 2, that
there exists a matrix $\mathbf{G}^{h}\in H^{1,h}(\Omega ^{h})$ such that $%
\widetilde{\mathbf{G}}^{h}\mid _{\partial \Omega ^{h}}=\mathbf{F}^{h}.$ In
addition, we assume that $\mathbf{G}^{h}$ is given with a noise of the level 
$\delta \in \left( 0,1\right) ,$ i.e.%
\begin{equation}
\Vert \mathbf{G}^{\ast ,h}-\mathbf{G}^{h}\Vert _{H^{1,h}(\Omega ^{h})}\leq
\delta .  \label{5.10}
\end{equation}%
Our main goal now is to estimate the difference between $\mathbf{U}_{\min
,\epsilon }^{h}$ and $\mathbf{U}^{\ast h}$ via $\delta $ and $\epsilon .$

\begin{lemma}
There exists a number $C_{2}=C_{2}(N,d,R,h_{1},a,b,L^{h})>0$ 
and a sufficiently small number $h_{y}^{0}=h_{y}^{0}(N,d,R,h_{1},a,b,L^{h})%
\in \left( 0,1\right) ,$ both depending only on listed parameters,
such that for $h_{x}\in \left[ h_{1},1\right) ,h_{y}\in \left( 0,h_{y}^{0}%
\right] $\ the following estimate is valid
\begin{equation}
\left\Vert L^{h}\left( \mathbf{Q}^{h}\right) \right\Vert _{L^{2,h}\left(
\Omega ^{h}\right) }^{2}\geq C_{2}\left\Vert \mathbf{Q}^{h}\right\Vert
_{L^{2,h}\left( \Omega ^{h}\right) }^{2}, \quad \mbox{for all } \widetilde{\mathbf{Q}}%
^{h}\in H_{0}^{1,h}(\Omega ^{h}).  \label{5.11}
\end{equation}%
\label{lem3}
\end{lemma}
\textbf{Proof.} Below $C_{2}>0$ denotes different constants depending on the
above listed parameters. Using the definition of the operator $L^{h}$ in (%
\ref{4.6}) as well as (\ref{5.4}), (\ref{5.5}) and the Cauchy-Schwarz
inequality, we obtain

\begin{align*}
\begin{split}
&\left\Vert L^{h}\left( \mathbf{Q}^{h}\right) \right\Vert _{L^{2,h}\left(
\Omega ^{h}\right) }^{2} = \left\Vert L^{h}\left( \mathbf{Q}^{h}\right)
e^{\lambda y}e^{-\lambda y}\right\Vert _{L^{2,h}\left( \Omega ^{h}\right)
}^{2}\geq \\ &e^{-2\lambda b}\left\Vert L^{h}\left( \mathbf{Q}^{h}\right)
e^{\lambda y}\right\Vert _{L^{2,h}\left( \Omega ^{h}\right) }^{2} \geq \frac{1}{2}e^{-2\lambda b}\left[ \left\Vert \left( \mathbf{Q}%
^{h}\right) _{y}^{\prime }e^{\lambda y}\right\Vert _{L^{2,h}\left( \Omega
^{h}\right) }^{2}-C_{2}\left\Vert \mathbf{Q}^{h}\right\Vert _{L^{2,h}\left(
\Omega ^{h}\right) }^{2}\right] .
\end{split}%
\end{align*}%
Choose $h_{y}^{0}\in \left( 0,1\right) $ so small that $1/\left(
h_{y}^{0}\right) ^{2}>64C_{2}$ and let $h_{y}\in \left( 0,h_{y}^{0}\right) .$
Set $\lambda =1/\left( 2h_{y}\right) .$ Then $\lambda h_{y}<1/2<1$ and also $%
\lambda ^{2}/8=1/\left( 32\left( h_{y}\right) ^{2}\right) >2C_{2}.$ Hence, by (\ref{5.8}) and (\ref{5.11}) it follows from the above inequality%
\begin{align*}
\begin{split}
\left\Vert L^{h}\left( \mathbf{Q}^{h}\right) \right\Vert _{L^{2,h}\left(
\Omega ^{h}\right) }^{2} &\geq \frac{1}{2}e^{-2\lambda b}\left[ \frac{\lambda
^{2}}{8}\left\Vert \mathbf{Q}^{h}e^{\lambda y}\right\Vert _{L^{2,h}\left( \Omega
^{h}\right) }^{2}-C_{2}\left\Vert \mathbf{Q}^{h}e^{\lambda y}\right\Vert
_{L^{2,h}\left( \Omega ^{h}\right) }^{2}\right] 
 \\ 
 & \geq \frac{1}{2} e^{-2\lambda b}\left( 2C_{2}\left\Vert \mathbf{Q}%
^{h}e^{\lambda y}\right\Vert _{L^{2,h}\left( \Omega ^{h}\right)
}^{2}-C_{2}\left\Vert \mathbf{Q}^{h}e^{\lambda y}\right\Vert _{L^{2,h}\left(
\Omega ^{h}\right) }^{2}\right) 
 \\
& =\frac{1}{2} e^{-2\lambda b}C_{2}\left\Vert \mathbf{Q}^{h}e^{\lambda
y}\right\Vert _{L^{2,h}\left( \Omega ^{h}\right) }^{2}.
\end{split}%
\end{align*}%
This estimate immediately implies (\ref{5.10}) with a new constant $C_{2}>0.$
$\square $ \newline

\begin{theorem}[Convergence rate of regularized solutions]
Assume that conditions of Theorem \ref{thm2} as well as %
(\ref{5.10}) and (\ref{5.11}) hold. 
Let $U_{\min, \epsilon }^{h}\in H^{1,h}(\Omega ^{h})$ be the unique minimizer of
the functional (\ref{4.9}) that satisfies boundary condition (%
\ref{4.7}) (see Theorem \ref{thm2}). Suppose that $h_{x}\in \lbrack h_{1},1)$\textit{\
and }$h_{y}\in \left( 0,h_{y}^{0}\right] ,$\textit{\ where the number }$%
h_{y}^{0}$\textit{\ is defined in Lemma \ref{lem3}. Then for any }$\epsilon >0$%
\textit{\ the following convergence rate of regularized solutions holds} 
\begin{equation}
\Vert \mathbf{U}_{\min ,\epsilon }^{h}-\mathbf{U}^{\ast h}\Vert
_{L^{2,h}(\Omega ^{h})}\leq C_{2}(\delta +\sqrt{\epsilon }\Vert \mathbf{U}%
^{\ast h}\Vert _{H^{1,h}(\Omega ^{h})}).  \label{5.12}
\end{equation}
\end{theorem}

\textbf{Proof.} Define the matrix $\widetilde{\mathbf{W}}_{\min ,\epsilon
}^{h}\in H_{0}^{1,h}(\Omega ^{h})$ as in Theorem \ref{thm2}, i.e. $%
\widetilde{\mathbf{W}}_{\min ,\epsilon }^{h}=\widetilde{\mathbf{U}}_{\min
,\epsilon }^{h}-\widetilde{\mathbf{G}}^{h}.$ Similarly define $\widetilde{%
\mathbf{W}}^{\ast ,h}=\widetilde{\mathbf{U}}^{\ast ,h}-\widetilde{\mathbf{G}}%
^{\ast ,h}\in H_{0}^{1,h}(\Omega ^{h}).$ Then (\ref{5.3}) is valid for $%
\mathbf{W}_{\min ,\epsilon }^{h}.$ As to $\mathbf{W}^{\ast ,h},$ (\ref{4.6})
and (\ref{4.7}) imply that the following analog of (\ref{5.3}) is valid for
all $\widetilde{\mathbf{P}}^{h}\in H_{0}^{1,h}(\Omega ^{h}):$%
\begin{equation*}
\left( L^{h}\left( \mathbf{W}^{\ast ,h}\right) ,L^{h}\left( \mathbf{P}%
^{h}\right) \right) +\epsilon \left[ \mathbf{W}^{\ast ,h},\mathbf{P}^{h}%
\right]  =-\left( L^{h}\left( \mathbf{G}^{\ast h}\right) ,L^{h}\left( 
\mathbf{P}^{h}\right) \right) 
 +\epsilon \left[ \mathbf{W}^{\ast ,h},\mathbf{P}^{h}\right] .
\end{equation*}%
Denote $\widetilde{\mathbf{V}}^{h}=\widetilde{\mathbf{W}}_{\min ,\epsilon
}^{h}-\mathbf{W}^{\ast ,h}\in H_{0}^{1,h}(\Omega ^{h}),\mathbf{Y}^{h}=%
\mathbf{G}^{h}-\mathbf{G}^{\ast h}$ and subtract (\ref{5.14}) from (\ref{5.3}%
). We obtain for all $\widetilde{\mathbf{P}}^{h}\in H_{0}^{1,h}(\Omega ^{h})$%
\begin{equation}
\left( L^{h}\left( \mathbf{V}^{h}\right) ,L^{h}\left( \mathbf{P}^{h}\right)
\right) +\epsilon \left[ \mathbf{V}^{h},\mathbf{P}^{h}\right] =-\left(
L^{h}\left( \mathbf{Y}^{h}\right) ,L^{h}\left( \mathbf{P}^{h}\right) \right)
-\epsilon \left[ \mathbf{W}^{\ast ,h},\mathbf{P}^{h}\right] .  \label{5.13}
\end{equation}%
Set in (\ref{5.13}) $\mathbf{P}^{h}=\mathbf{V}^{h}$ and use the
Cauchy-Schwarz inequality. We obtain 
\begin{equation*}
\begin{split}
\left\Vert L^{h}\left( \mathbf{V}^{h}\right) \right\Vert _{L^{2,h}\left(
\Omega ^{h}\right) }^{2}& \leq \left\Vert L^{h}\left( \mathbf{Y}^{h}\right)
\right\Vert _{L^{2,h}\left( \Omega ^{h}\right) }^{2}+\epsilon \left\Vert 
\mathbf{W}^{\ast ,h}\right\Vert _{H^{1,h}\left( \Omega ^{h}\right) }^{2}
\\
& \leq C_{2}\delta ^{2}+\epsilon \left\Vert \mathbf{W}^{\ast ,h}\right\Vert
_{H^{1,h}\left( \Omega ^{h}\right) }^{2}.
\end{split}%
\end{equation*}%
Next, using (\ref{5.11}) and (\ref{5.13}), we obtain 
\begin{equation}
\left\Vert \mathbf{V}^{h}\right\Vert _{L^{2,h}\left( \Omega ^{h}\right)
}^{2}\leq C_{2}\left( \delta ^{2}+\epsilon \left\Vert \mathbf{W}^{\ast
,h}\right\Vert _{H^{1,h}\left( \Omega ^{h}\right) }^{2}\right) .  \label{5.14}
\end{equation}%
The target estimate (\ref{5.12}) follows immediately from (\ref{5.14}). $%
\square $ \newline

\begin{theorem}[Lipschitz stability and uniqueness.]
Suppose that there exist two matrices $\mathbf{G}_{1}^{h},\mathbf{G}%
_{2}^{h}\in H^{1,h}(\Omega ^{h})$ such that $\widetilde{\mathbf{G}}%
_{1}^{h}\mid _{\partial \Omega ^{h}}=\mathbf{F}_{1}^{h}$ \textit{and }$%
\widetilde{\mathbf{G}}_{2}^{h}\mid _{\partial \Omega ^{h}}=\mathbf{F}%
_{2}^{h},$ \textit{where} $\mathbf{F}_{1}^{h}$ \textit{and} $\mathbf{F}_{2}^{h}$ 
are two different boundary conditions in (\ref{4.7}).
Suppose that there exist solutions $\widetilde{\mathbf{U}}_{1}^{h}\in
H^{1,h}\left( \Omega ^{h}\right) $ and $\widetilde{\mathbf{U}}_{2}^{h}\in
H^{1,h}\left( \Omega ^{h}\right) $ \textit{of boundary value problem %
\normalfont{(\ref{4.6})-(\ref{4.7})} with boundary conditions }$\mathbf{F}%
_{1}^{h}$ \textit{and} $\mathbf{F}_{2}^{h}$ \textit{respectively. Assume that }$%
h_{x}\in \lbrack h_{1},1)$\textit{\ and }$h_{y}\in \left( 0,h_{y}^{0}\right] ,$%
\textit{\ where the number }$h_{y}^{0}$\textit{\ is defined in Lemma \ref{lem3}.
Then the following Lipschitz stability estimate is valid}%
\begin{equation}
\left\Vert \mathbf{U}_{1}^{h}-\mathbf{U}_{2}^{h}\right\Vert _{L^{2,h}\left(
\Omega ^{h}\right) }\leq C_{2}\left\Vert \mathbf{G}_{1}^{h}-\mathbf{G}%
_{2}^{h}\right\Vert _{L^{2,h}\left( \Omega ^{h}\right) }.  \label{5.15}
\end{equation}%
\textit{Next, suppose that} $\mathbf{F}_{1}^{h}=\mathbf{F}_{2}^{h},$ \textit{but
the existence of the function} $\widetilde{\mathbf{G}}^{h}$ \textit{is not
assumed. Then} $\mathbf{U}_{1}^{h}=\mathbf{U}_{2}^{h},$ \textit{where} $%
\mathbf{U}_{1}^{h}$,$\mathbf{U}_{2}^{h}\in H^{1,h}\left( \Omega ^{h}\right) $
\textit{are two possible solution of boundary value problem (\ref{4.6}), (\ref%
{4.7}).}
\end{theorem}

\textbf{Proof}. Since $\widetilde{\mathbf{U}}_{1}^{h}$ and $\widetilde{%
\mathbf{U}}_{2}^{h}$ are two exact solutions of problem (\ref{4.6}), (\ref%
{4.7}) with two different boundary conditions, then by (\ref{5.12})%
\begin{equation}
\left( L^{h}\left( \mathbf{U}_{i}^{h}\right) ,L^{h}\left( \mathbf{P}%
^{h}\right) \right) =-\left( L^{h}\left( \mathbf{G}_{i}^{h}\right)
,L^{h}\left( \mathbf{P}^{h}\right) \right) ,\forall \mathbf{P}^{h}\in
H_{0}^{1,h}\left( \Omega ^{h}\right) ,i=1,2,  \label{5.16}
\end{equation}%
where $\widetilde{\mathbf{W}}_{i}^{h}=\widetilde{\mathbf{U}}_{i}^{h}-%
\widetilde{\mathbf{G}}_{i}^{h}.$\ Setting $\mathbf{S}^{h}=\mathbf{W}_{1}^{h}-%
\mathbf{W}_{2}^{h},\mathbf{X}^{h}=\mathbf{G}_{1}^{h}-\mathbf{G}_{2}^{h}$ and
then setting $\mathbf{P}^{h}=\mathbf{S}^{h},$ we obtain from (\ref{5.16})
and (\ref{4.8}) 
\begin{equation*}
\left\Vert L^{h}\left( \mathbf{S}^{h}\right) \right\Vert _{L^{2,h}\left(
\Omega ^{h}\right) }^{2}\leq \left\Vert L^{h}\left( \mathbf{X}^{h}\right)
\right\Vert _{L^{2,h}\left( \Omega ^{h}\right) }^{2}\leq C_{2}\left\Vert 
\mathbf{X}^{h}\right\Vert _{L^{2,h}\left( \Omega ^{h}\right) }^{2}.
\end{equation*}%
Hence, by (\ref{5.11}) $\left\Vert \mathbf{S}^{h}\right\Vert _{L^{2,h}\left(
\Omega ^{h}\right) }^{2}\leq C_{2}\left\Vert \mathbf{X}^{h}\right\Vert
_{L^{2,h}\left( \Omega ^{h}\right) }^{2}.$ Therefore,%
\begin{equation}
\left\Vert \left( \mathbf{U}_{1}^{h}-\mathbf{U}_{2}^{h}\right) -\left( 
\mathbf{G}_{1}^{h}-\mathbf{G}_{2}^{h}\right) \right\Vert _{L^{2,h}\left(
\Omega ^{h}\right) }\leq C_{2}\left\Vert \mathbf{G}_{1}^{h}-\mathbf{G}%
_{2}^{h}\right\Vert _{L^{2,h}\left( \Omega ^{h}\right) }.  \label{5.17}
\end{equation}%
Thus, (\ref{5.15}) follows from (\ref{5.17}) and the triangle inequality. As
to the uniqueness part, since $\mathbf{F}_{1}^{h}=\mathbf{F}_{2}^{h},$ then
we extend the boundary condition $\left( \mathbf{F}_{1}^{h}-\mathbf{F}%
_{2}^{h}\right) =0$ in the domain $\Omega ^{h}$ as $\mathbf{G}_{1}^{h}-%
\mathbf{G}_{2}^{h}\equiv 0.$ Hence, (\ref{5.15}) implies that $\mathbf{U}%
_{1}^{h}-\mathbf{U}_{2}^{h}\equiv 0.$ \ $\square $

\begin{remark}
Due to the ill-posedness of Problem \ref{ip}, we cannot prove convergence of our solutions to the correct one as $%
N \to \infty $. We note that the truncated Fourier series is used
both quite often and quite successfully in numerical methods for many
inverse problems. Although convergences at $N \to \infty $ are
not proven in many cases, numerical results are usually good ones,
see, e.g. \cite{N} for the attenuated tomography with complete data, \cite%
{Kab05} for the 2D version of the Gelfand-Levitan method,
\cite{LiNguyen:arxiv2019} for the inverse problem  of finding initial condition of heat equation, \cite{NguyenLiKlibanov:arxiv2019} for the inverse source problem for the Helmholtz equation, and \cite%
{KTsiam15,KK18,KLZ18} for the convexification.
\end{remark}

\section{Numerical Implementation}

\label{sec 6}

In this section, we describe the numerical implementation of the
minimization procedure for the functional $J_{\epsilon }^{h}$. While inverting the matrix $M_N$ of (\ref{3.6}) is convenient for the convergence analysis, we have discovered that it is better in real computations not to invert while still considering a problem which is equivalent to the problem (\ref{4.6}),(\ref{4.7}). Thus, we we consider the functional (\ref{4.9}) in a slightly different form: 
\begin{equation*}
\begin{split}
J_{\epsilon _{1},\epsilon _{2}}^{h}(\normalfont{\widetilde{\mathbf{U}}^h}%
)=\Vert \normalfont{(M_N-\mathbf{A}^h) ({\mathbf{U}}^h)^{\prime}_y}%
&-\mathbf{B}^{h} \normalfont{({\mathbf{U}}^h)^{\prime}_x} -\mathbf{C}^{h} \normalfont{{%
\mathbf{U}}^h}\Vert _{L^{2,h}({\Omega }^{h})}^{2}+ \\
& +\epsilon _{1}\Vert \mathbf{U}^{h}\Vert _{L^{2,h}({\Omega }%
^{h})}^{2}+\epsilon _{2}\Vert \nabla \mathbf{U}^{h}\Vert _{L^{2,h}({\Omega }^{h})}^{2},
\end{split}%
\end{equation*}%
where $\mathbf{A}^{h},\mathbf{B}^{h},\mathbf{C}^{h}$ are operators (\ref{3.7})-(\ref{3.9}), with the domain $\overline{\Omega}_h$. Moreover, in contrast to the original functional, we use in our computations two
regularization parameters $\epsilon _{1}$ and $\epsilon _{2}$, instead of
just one parameter $\epsilon $. This yields better reconstruction results.
The regularization parameters $\epsilon _{1},\epsilon _{2}$ in our numerical
tests were found by a trial and error procedure. They were the same for all the tests we have conducted.

To minimize the functional $J_{\epsilon _{1},\epsilon _{2}}^{h}(%
\normalfont{\widetilde{\mathbf{U}}^h}),$ we first have to simulate the
boundary data via solving Problem \ref{pro forward}. 
We discuss the
solution of this forward problem in Appendix. Using this solution, we
generate the noisy data as, see (\ref{2.8}), (\ref{2.9}) and (\ref{5.10}):%
\begin{equation*}
F\left( \mathbf{x},\alpha \right) =\left\{ 
\begin{array}{c}
u_{\delta }^{comp}(\mathbf{x},\alpha )=u^{comp}(\mathbf{x},\alpha )(1+\delta
(2\mathrm{rand}(\mathbf{x})-1)),\text{ }\mathbf{x}\in \partial \Omega
\setminus \partial \Omega ^{+}, \\ 
0,\text{ }\mathbf{x}\in \partial \Omega ^{+}%
\end{array}%
\right. .
\end{equation*}%
where $u^{comp}(\mathbf{x},\alpha )$ is the boundary data computed via the
solution of the forward problem, $\delta >0$ is the noise level and $\mathrm{%
rand(\cdot )}$ is the function that generates uniformly distributed random
numbers in the interval $[0,1]$. For example, $\delta =0.6$ corresponds to
the 60\% noise level in the data.

We use finite the difference approximations (\ref{4.4}), (\ref{4.5}) on the
grid with $h_{x}=h_{y},$ $M_{x}=M_{y}$. We rewrite the functional $%
J_{\epsilon _{1},\epsilon _{2}}^{h}$ in the following discrete form 
\begin{align*}
& J_{\epsilon _{1},\epsilon _{2}}^{h} (\normalfont{\widetilde{\mathbf{U}}^h}) = \\
&\displaystyle h_{\mathbf{x}}^{2}\sum_{i,j=1}^{M_{x}-1}\left( (M_N-\mathbf{A}^h)_{ij} \frac{\mathbf{U}%
_{i,j+1}^{h}-\mathbf{U}_{i,j}^{h}}{h_{x}}-(\mathbf{B}^{h})_{ij}\frac{%
\mathbf{U}_{i+1,j}^{h}-\mathbf{U}_{i,j}^{h}}{h_{x}}-(\mathbf{C}^{h})_{ij}%
\mathbf{U}_{i,j}^{h}\right) ^{2}
\\
&+ \displaystyle\epsilon _{1}\displaystyle h_{x}^{2}\sum_{i,j=1}^{M_{x}-1}%
\left( \mathbf{U}_{i,j}^{h}\right) ^{2}+\displaystyle\epsilon _{2}%
\displaystyle h_{x}^{2}\sum_{i,j=1}^{M_{x}-1}\left( \frac{|\mathbf{U}%
_{i,j+1}^{h}-\mathbf{U}_{i,j}^{h}|^{2}}{h_{x}}+\frac{|\mathbf{U}_{i+1,j}^{h}-%
\mathbf{U}_{i,j}^{h}|^{2}}{h_{x}}\right) .
\end{align*}%
Denote $u_{m}(x_{i},y_{j})=u_{m}^{i,j}$. Since $\mathbf{U}
:=(u_{1}(x,y),\dots ,u_{N}(x,y))^{T},$ then%
\begin{align*}
& J_{\epsilon _{1},\epsilon _{2}}^{h} (\widetilde{\mathbf{U}}^{h})
= \\ &h_{x}^{2}\sum_{i,j=1}^{M_{x}-1}\sum_{m=1}^{N}\left((M_N-\mathbf{A}^h)_{ij} \frac{%
u_{m}^{i,j+1}-u_{m}^{i,j}}{h_{x}}-(\mathbf{B}^{h})_{ij}\frac{%
u_{m}^{i+1,j}-u_{m}^{i,j}}{h_{x}}-(\mathbf{C}^{h})_{ij}u_{m}^{i,j}%
\right)^{2} \\
&+ \displaystyle\epsilon_{1}\displaystyle h_{x}^{2}\sum_{i,j=1}^{M_{x}-1}%
\sum_{m=1}^{N}\left( u_{m}^{i,j})\right) ^{2}+\displaystyle\epsilon _{2}%
\displaystyle h_{x}^{2}\sum_{i,j=1}^{M_{x}-1}\sum_{m=1}^{N}\left( \frac{%
|u_{m}^{i,j+1}-u_{m}^{i,j}|^{2}}{h_{x}}+\frac{|u_{m}^{i+1,j}-u_{m}^{i,j}|^{2}%
}{h_{x}}\right) .
\end{align*}%
Introduce the \textquotedblleft lined-up" versions of the matrices $%
\normalfont{\widetilde{\mathbf{U}}^h},M_N-\mathbf{A}^h, \mathbf{B}^{h},\mathbf{C}^{h}$. The \newline $(M_{x}+1)^{2}N$ dimensional vector $\mathcal{U},$ 
\begin{equation}
\mathcal{U}_{\mathfrak{m}}=u_{m}(x_{i},y_{j})\quad 1\leq i,j\leq
M_{x}+1,1\leq m\leq N,  \label{Ulineup}
\end{equation}%
and the $(M_{x}+1)^{2}N \times (M_{x}+1)^{2}N$ dimensional matrices $\mathcal{A}^{h},\mathcal{B}^{h}$, $\mathcal{C}^h$, corresponding to $M_N-\mathbf{A}^h, \mathbf{B}^{h},\mathbf{C}^{h}$, where 
\begin{equation}
\mathfrak{m}=(i-1)(M_{x}+1)N+(j-1)N+m.  
\label{line up index}
\end{equation}%
We introduce the map 
\begin{equation*}
\{1,\dots ,M_{x}+1\}\times \{1,\dots ,M_{x}+1\}\times \{1,\dots
,N\}\rightarrow \{1,\dots ,(M_{x}+1)^{2}N\}
\end{equation*}%
that sends $(i,j,m)$ to $\mathfrak{m}$ as in \eqref{line up index} is onto
and one-to-one. The functional $J_{\epsilon _{1},\epsilon _{2}}^{h}(%
\normalfont{\widetilde{\mathbf{U}}^h})$ is rewritten in terms of the
lined-up vector $\mathcal{U}$ \ as 
\begin{equation}
\mathcal{J}_{\epsilon _{1},\epsilon _{2}}^{h}(\mathcal{U})=h_{x}^{2}\left( |%
\mathcal{L}\mathcal{U}|^{2}+\epsilon _{1}|\mathcal{U}|^{2}+\epsilon _{2}|%
\mathcal{D}_{x}\mathcal{U}|^{2}+\epsilon _{2}|\mathcal{D}_{x}\mathcal{U}%
|^{2}\right) ,  \label{mathcalJ}
\end{equation}%
where $\mathcal{D}_{x}$ and $\mathcal{D}_{y}$ are the matrices that provide
the finite difference analogs of the partial derivatives of $\mathcal{U}$
with respect to $x$ and $y$, defined similarly to (\ref{4.4}),(\ref{4.5}). $%
\mathcal{L}$ is the $(M_{x}+1)^{2}N\times (M_{x}+1)^{2}N$ matrix defined as
follows. For each 
\begin{equation}
\mathfrak{m}=(i-1)(M_{x}+1)N+(j-1)N+m,2\leq i,j\leq M_{x},1\leq m\leq N
\label{index}
\end{equation}%

\begin{enumerate}
\item $\mathcal{L}_{\mathfrak{m}\mathfrak{n}}=(-(\mathcal{A}^h)_{%
\mathfrak{m}\mathfrak{n}}+(\mathcal{B}^h)_{%
\mathfrak{m}\mathfrak{n}})/h_x - (\mathcal{C}^h)_{\mathfrak{m}\mathfrak{n}%
},$ if $\mathfrak{m}$ corresponds to $(i, j, n)$ in the sense of %
\eqref{index} for each $n\in \{1,\dots ,N\}$,

\item $\mathcal{L}_{\mathfrak{m}\mathfrak{n}}=(\mathcal{A}^h)_{%
\mathfrak{m}\mathfrak{n}}/h_x,$ if $\mathfrak{m}$
corresponds to $(i, j+1, n)$ in the sense of \eqref{index} for each $n\in \{1,\dots ,N\}$,

\item $\mathcal{L}_{\mathfrak{m}\mathfrak{n}}=-(\mathcal{B}^h)_{\mathfrak{m%
}\mathfrak{n}}/h_x,$ if $\mathfrak{m}$ corresponds to $(i+1, j, n)$ in the
sense of \eqref{index} for each $n\in \{1,\dots ,N\}$,

\item $\mathcal{L}_{\mathfrak{m}\mathfrak{n}}=0$ otherwise.
\end{enumerate}
\vspace{0.1in}

Next, we consider the \textquotedblleft lined-up" version of the boundary
condition (\ref{4.7}). Let $\mathcal{D}$ be the $(M_{x}+1)^{2}N\times
(M_{x}+1)^{2}N$ diagonal matrix with $\mathfrak{m}^{\mathrm{th}}$ diagonal
entries taking value $1$ while the others equal $0$. This Dirichlet boundary
constraint of the vector $\mathcal{U}$ becomes $\mathcal{D}\mathcal{U}=%
\tilde{\mathcal{F}}.$ Here, the vector $\tilde{\mathcal{F}}$ is the
\textquotedblleft lined-up" vector of the data $F_{N}$ in the same manner
when we defined $\mathcal{U}$, see (\ref{Ulineup}). 
We solve the Inverse Source Problem by computing the vector $\mathcal{U}$, subject to constraint $\mathcal{D}\mathcal{U}=\tilde{\mathcal{F}}$, such that 
\begin{equation}
\mathcal{L}_{\mu }\mathcal{U}=\left( \mathcal{L}^{T}\mathcal{L}+\epsilon _{1}%
\text{Id}+\epsilon _{2}\mathcal{D}_{x}^{T}\mathcal{D}_{x}+\epsilon _{2}%
\mathcal{D}_{y}^{T}\mathcal{D}_{y}\right) \mathcal{U}=\vec{0},
\label{mathcalU}
\end{equation}%
which is equivalent to the minimization of the functional \eqref{mathcalJ}.

The knowledge of $\mathcal{U}$ yields the knowledge of $\normalfont{%
\widetilde{\mathbf{U}}^h}$ via \eqref{Ulineup}. Denote the result obtained
by the procedure of this section as $\normalfont{\widetilde{%
\mathbf{U}}^{comp}}=(u_{1}^{comp}(x,y),\dots ,u_{N}^{comp}(x,y))^{T}$. Using
this vector function, we calculate function $u^{comp}(x,y,\alpha )$
according via (\ref{3.1}). Next, the reconstructed function $f^{comp}(x,y)$
is determined as the averaged over $\alpha $ value of the source function $%
f^{comp}(x,y,\alpha )$ calculated via the substitution of \newline
$u^{comp}(x,y,\alpha )$ in (\ref{2.3}).

These arguments lead to the Algorithm \ref{alg 1} for solving Problem \ref{ip} 

\begin{algorithm}
\caption{\label{alg 1}The procedure to solve Problem \ref{ip}}\label{euclid}
\begin{algorithmic}[1]
    \State\, Choose a number $N$. Construct the functions $\Psi_m$, $1 \leq m \leq N,$ in Section 3.1 and compute the matrix $M_N$ as in (\ref{matrix M}).
    \State\,  Calculate the boundary data $\mathbf{F}^h$ for the vector valued function $\widetilde{\mathbf{U}}^h$ on $\partial \Omega^h$ via solving Problem \ref{pro forward} with $f = f_{true}$.
    \State\, \label{step quasi} Find an approximate solution of \eqref{4.6}-\eqref{4.7} by the quasi-reversibility method.
    \State\, \label{step vcomp} Having $\widetilde{\mathbf{U}}^h$, calculate $u_{comp}(\mathbf{x}, \alpha)$ for $\mathbf{x} \in \Omega^h$ via \eqref{3.1}.
    \State\, \label{step fcomp} Compute the reconstructed function $f_{comp}$ by \eqref{2.3}.
\end{algorithmic}
\end{algorithm}

\subsection{Numerical tests}

We test our method using numerical simulations for different types of
absorption and scattering coefficients. Test 1 demonstrates the stability of
the solution in the \textquotedblleft no scattering" model, which
corresponds to the $\mu _{s}(\mathbf{x},\alpha )\equiv 0$ value of the
scattering coefficient. Moreover, Test 1 is used to find the optimal
parameters: $\epsilon _{1},\epsilon _{2},d,N,$ which we use in subsequent
tests. The optimal values for those listed parameters are: $\epsilon
_{1}=0.1,\epsilon _{2}=0.01,d=5,N=12$.

In Test 2 we consider the \textquotedblleft uniform scattering" version of
the RTE for the object with a non-smooth boundary, where $\mu _{s}(\mathbf{x}%
)\not\equiv 0,K(\mathbf{x},\alpha ,\beta )=1/(2d).$ 

Test 3 demonstrates the
performance of our method for both sophisticated form of the absorption coefficient $\mu _{a}(\mathbf{x})$ and \textquotedblleft strongly forward-peaked scattering" case, i.e. the case when the scattered particles move in the direction preferentially close to the one they were moving. That corresponds to the case when $\mu _{s}(\mathbf{x})\not\equiv 0$ and $K(\mathbf{x},\alpha ,\beta )\not=0$ whenever $\vert \alpha - \beta \vert \leq C$, for some constant $C > 0$. We choose our scattering phase function $K(\mathbf{x},\alpha,\beta)$ to be a 2-dimensional Henyey-Greenstein function, which, possesses the property mentioned above and is convenient to describe strongly forward-peaked scattering, according to \cite{Heino}.

All numerical tests were provided on the domain $\overline{\Omega }^{h}$ defined in Section 4 on the uniform $100\times 100$ grid (\ref{4.1}),(\ref{4.2}) with $a=1,b=3,d=5$, 
$R=1,$ see (\ref{2.1}), (\ref{2.2}). We use the uniform grid for $\alpha $
as well 
\begin{equation*}
-d=\alpha _{0}<\alpha _{1}<\dots <\alpha _{M_{\alpha }}=d,\hspace{0.5em}%
\text{where}\hspace{0.5em}M_{\alpha }=50.
\end{equation*}

Similarly with \cite{KN}, we apply a 2-step post-processing procedure. Define \newline $m = \max\displaylimits_{\mathbf{x} \in \Omega^h}(f^{comp}_{\delta}(\mathbf{x}))$. Then, the first step of the procedure removes undesirable artifacts by setting 
\begin{equation*}
\tilde{f}^{comp}_{\delta}(\mathbf{x})=\left\{ 
\begin{array}{rl}
f^{comp}_{\delta}(\mathbf{x}) & \mbox{if } f^{comp}_{\delta}(\mathbf{x}) > 0.2 m \\ 
0 & \mbox{otherwise,}%
\end{array}%
\right. 
\end{equation*}%
In the second step, we smooth the obtained $\tilde{f}^{comp}_{\delta}(\mathbf{x})$ function out. For every grid point $\mathbf{x}$ of $\Omega^h$, the value of $\tilde{f}^{comp}_{\delta}(\mathbf{x})$ at that point is replaced by the mean value $\hat{f}^{comp}_{\delta}(\mathbf{x})$ of the function over neighboring grid points. For brevity, we use $f^{comp}_{\delta}(\mathbf{x})$ notation instead of $\hat{f}^{comp}_{\delta}(\mathbf{x})$ below.
In Tests 1-3 we display the exact functions $%
f_{true}\left( x,y\right) $ and computed functions $f^{comp}(x,y),f_{\delta
}^{comp}(x,y)$ before and after the use of the post processing procedure.

\subsection{Numerical results}
\begin{enumerate}
\item \textbf{Test 1. Circle-shaped smooth inclusion, no scattering} \newline

The function $f_{\mathrm{true}}$ is a smoothed circle of the radius $r=0.3$,
centered at $(0,2)$, compactly supported in $\Omega ^{h}$. $f_{\mathrm{true}%
} $ is depicted on the Figure \ref{fig model 1}(a). 
\begin{equation*}
\mu _{s}(\mathbf{x})\equiv 0,\quad \mu _{a}(\mathbf{x})=\left\{ 
\begin{array}{rl}
0.1 & \mbox{if }(x^{2}+y^{2})<0.8, \\ 
0 & \mbox{otherwise,}%
\end{array}%
\right.
\end{equation*}

The numerical solution for this case is displayed in Figure \ref{fig model 1}.

\begin{figure}
\begin{center}
\subfloat[The function $f_{\rm true}$]{\includegraphics[width = .3\textwidth]{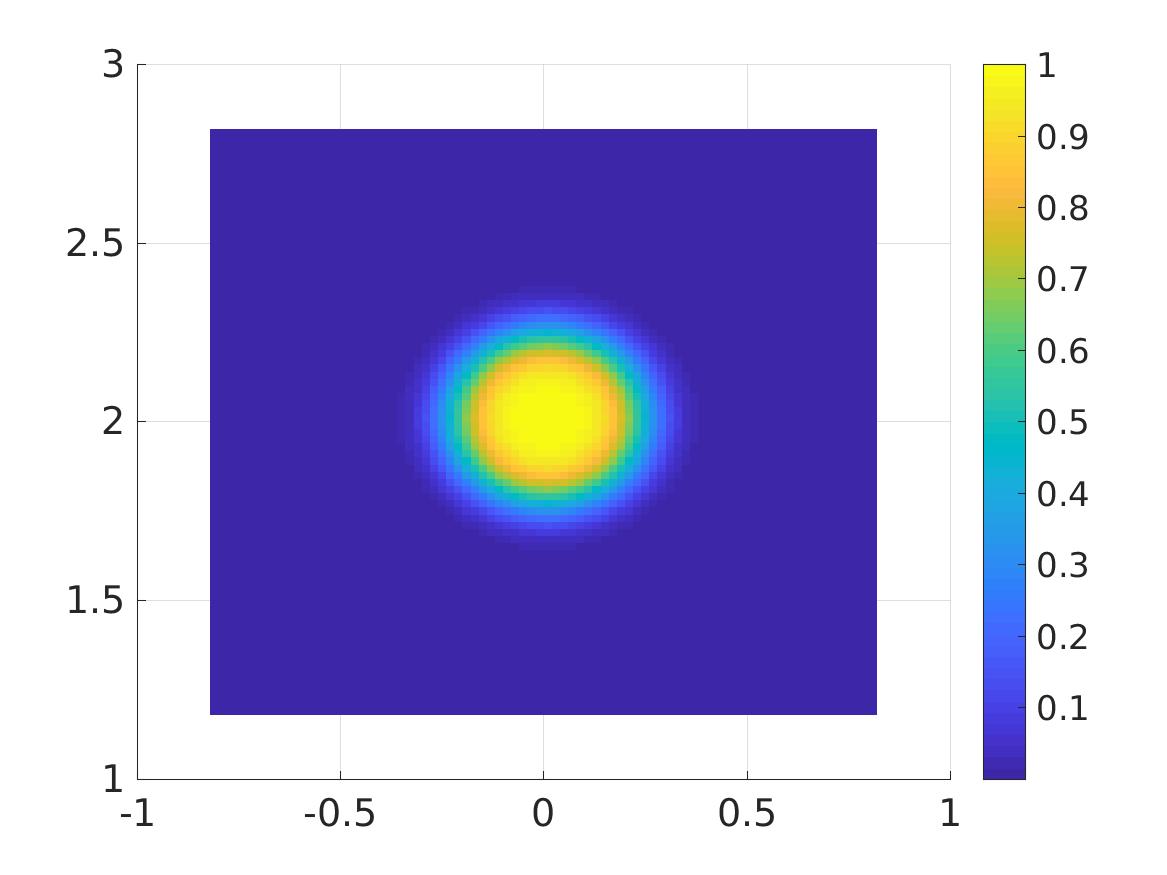}}\quad 
\subfloat[The reconstructed function $f^{\mathrm{comp}}$, no noise]{\includegraphics[width = .3\textwidth]{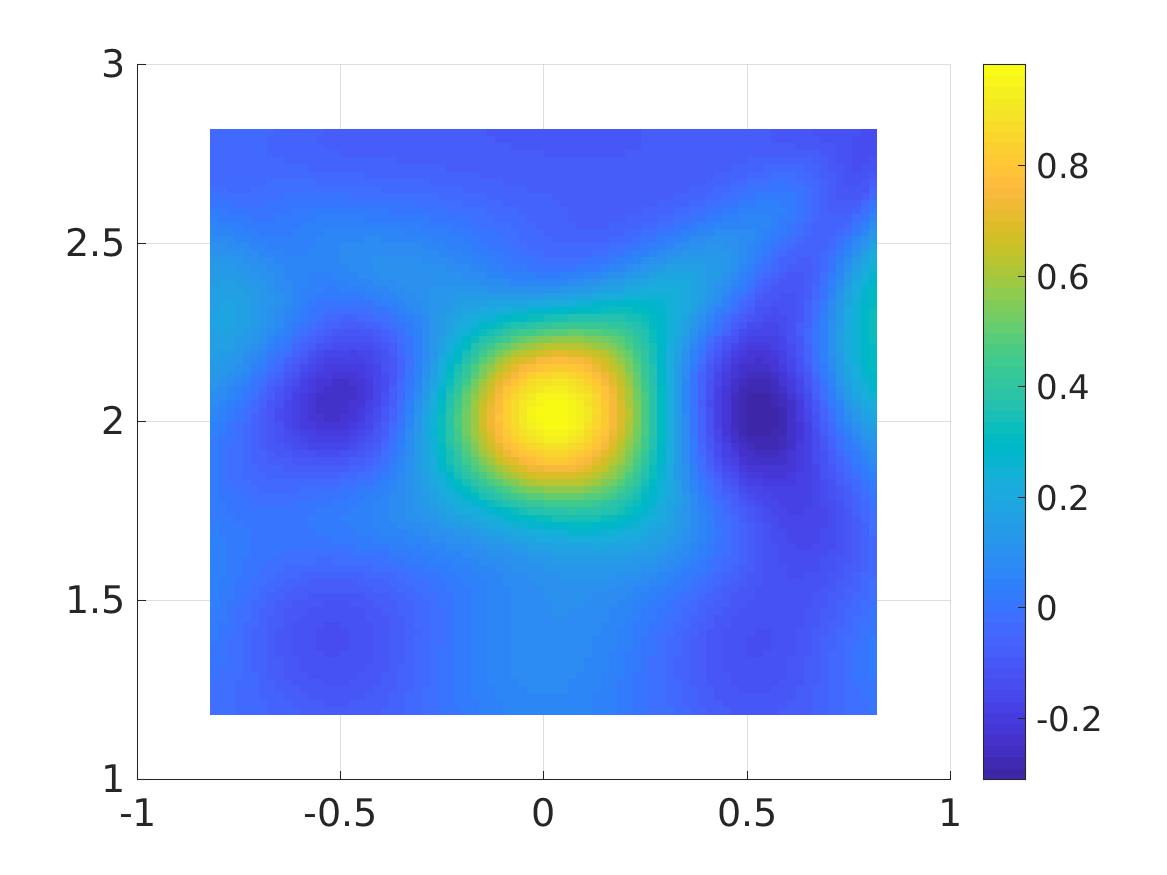}}  \quad
\subfloat[The reconstructed 
function $f^{\mathrm{comp}}_{\delta}$, noise 90\%]{\includegraphics[width = .3\textwidth]{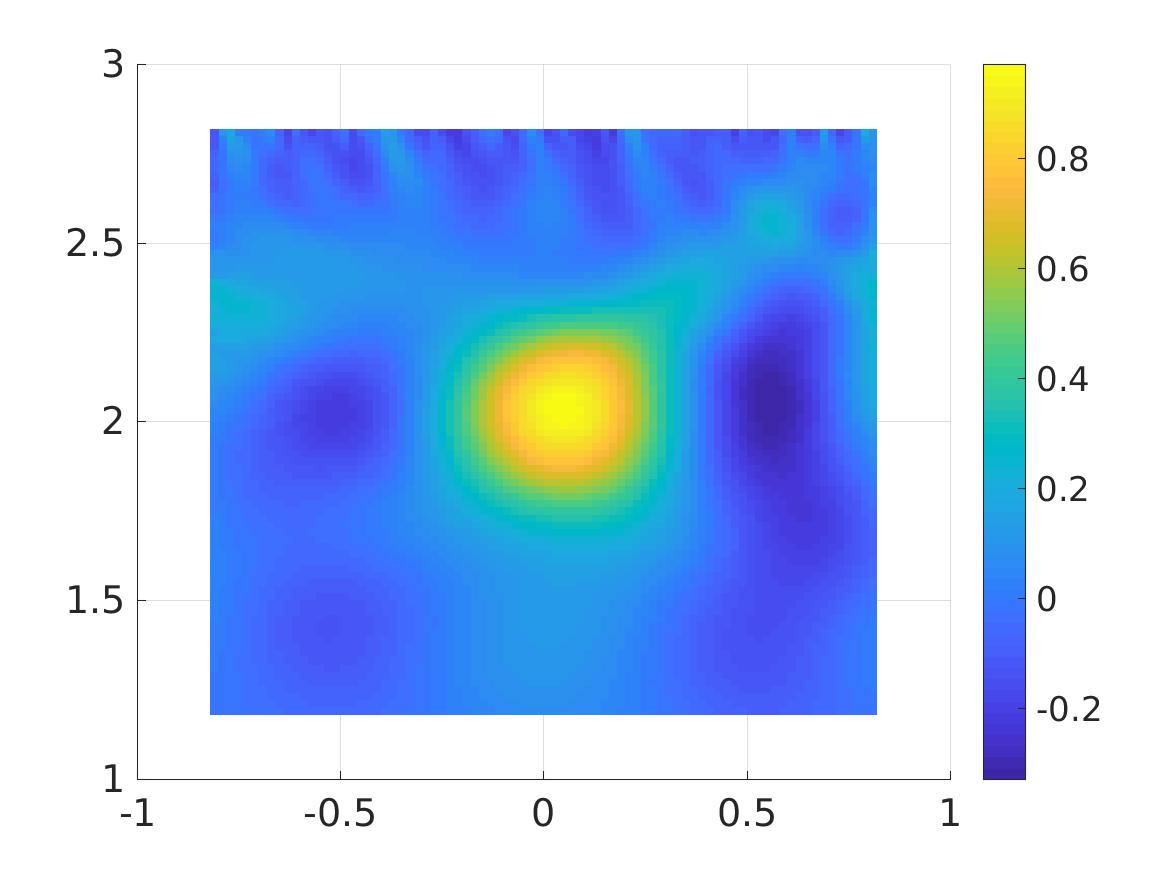}}
\par
\subfloat[\label{incomplete1} The incomplete tomographic data, no noise]{\includegraphics[width = .3\textwidth]{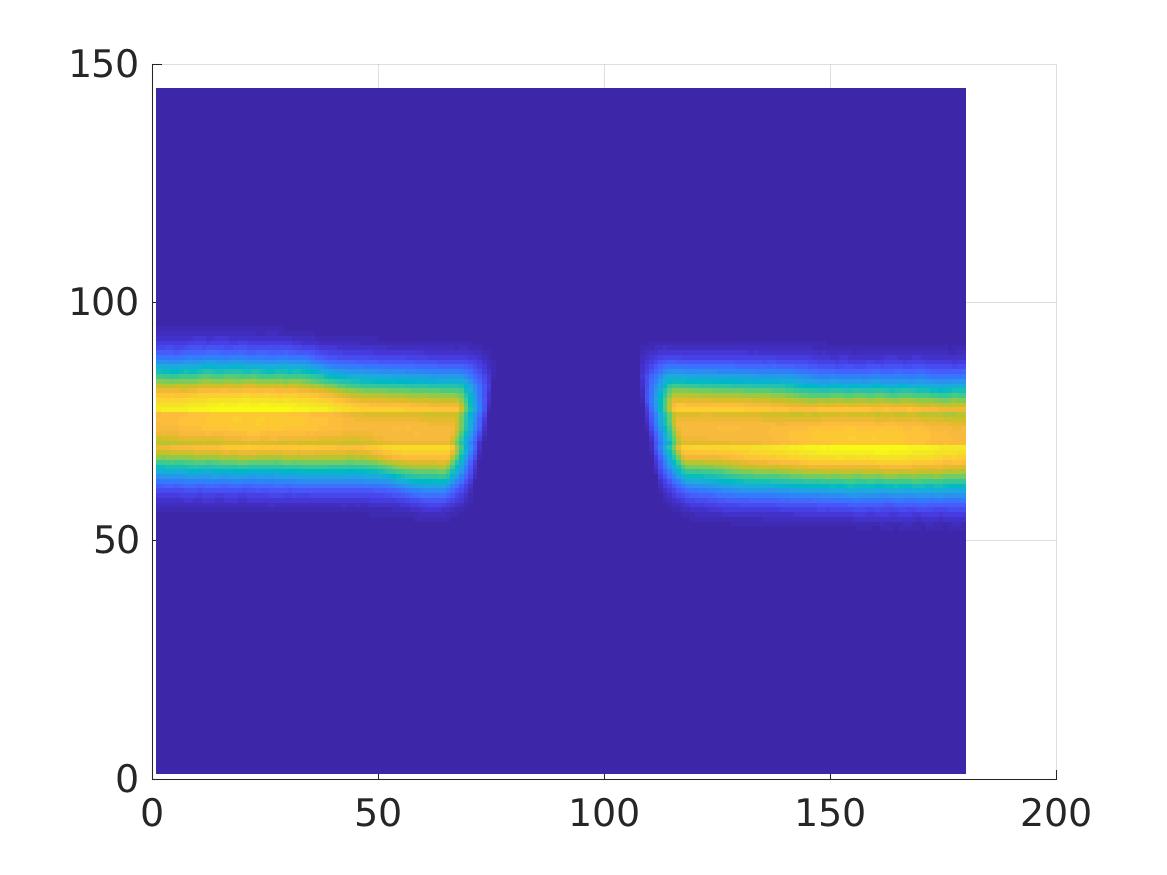}} \quad 
\subfloat[Post-processed $f^{\mathrm{comp}}$,  no noise]{\includegraphics[width = .3\textwidth]{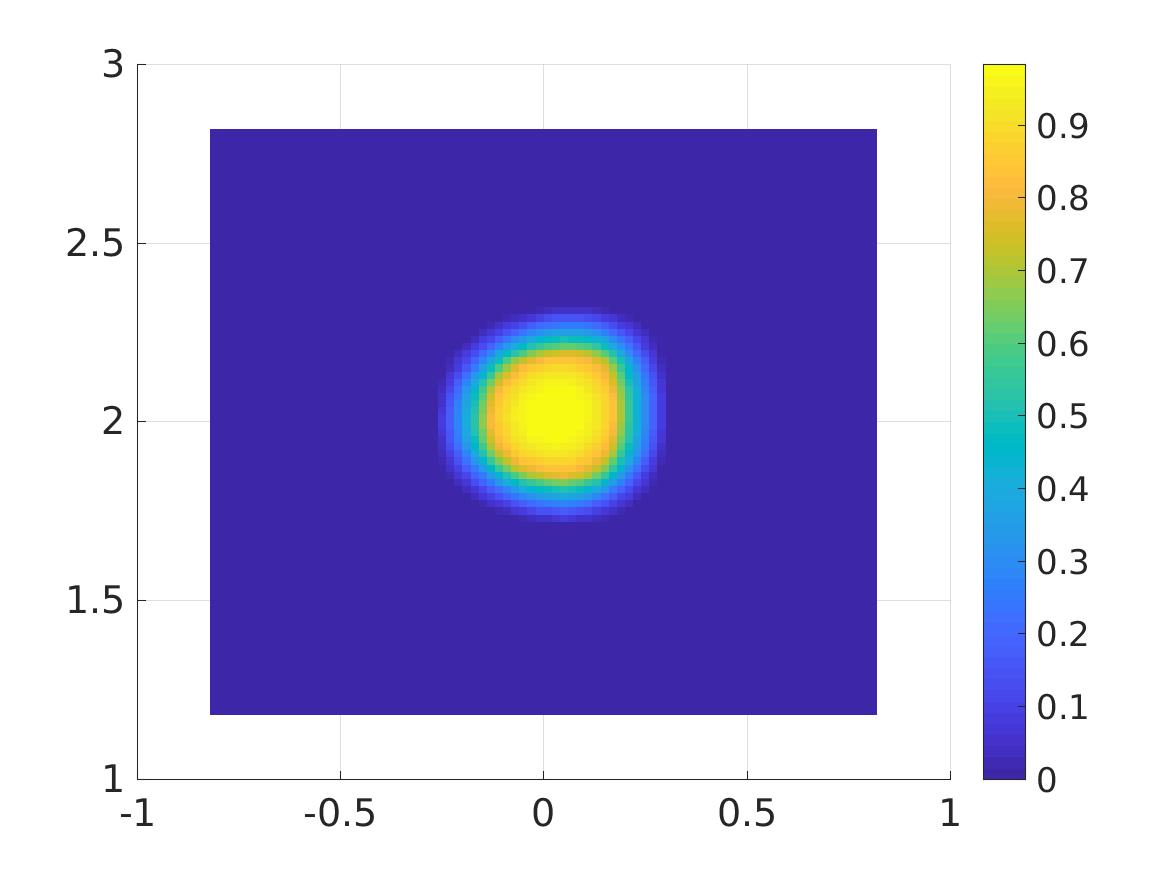}} \quad
\subfloat[Post-processed $f^{\mathrm{comp}}_{\delta}$, noise level 90\%]{\includegraphics[width = .3\textwidth]{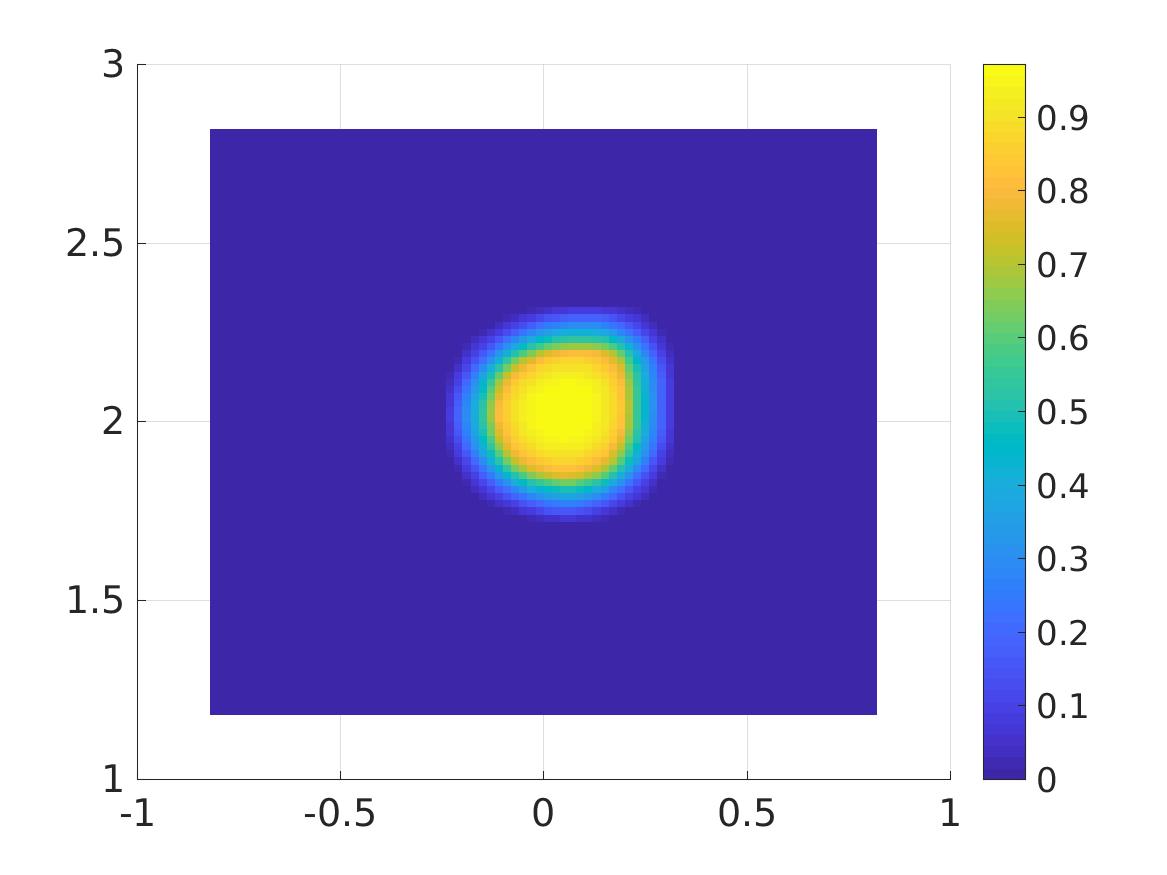}}
\end{center}
\caption{\textit{Test 1. The true and reconstructed source functions for
circle-shaped smooth inclusion without scattering.}}
\label{fig model 1}
\end{figure}

\item \textbf{Test 2. X-shaped non-smooth inclusion, uniform scattering} 

The function $f_{\mathrm{true}}$ is depicted on the Figure \ref{fig model 2}
(a). We define absorption and scattering coefficients as 
\begin{equation*}
\mu _{a}(\mathbf{x})=\left\{ 
\begin{array}{rl}
0.1 & \mbox{if }(x^{2}+y^{2})<0.8, \\ 
0 & \mbox{otherwise,}%
\end{array}%
\right. \mu _{s}(\mathbf{x},\alpha )=\left\{ 
\begin{array}{rl}
0.01 & \mbox{if }(x^{2}+y^{2})<0.8, \\ 
0 & \mbox{otherwise,}%
\end{array}%
\right.
\end{equation*}%
and the constant scattering phase function $K(\mathbf{x},\alpha ,\beta
)=1/(2d).$ The reconstruction is displayed in Figure \ref{fig model 2}. 
\begin{figure}
\begin{center}
\subfloat[The function $f_{\rm true}$]{\includegraphics[width =
.3\textwidth]{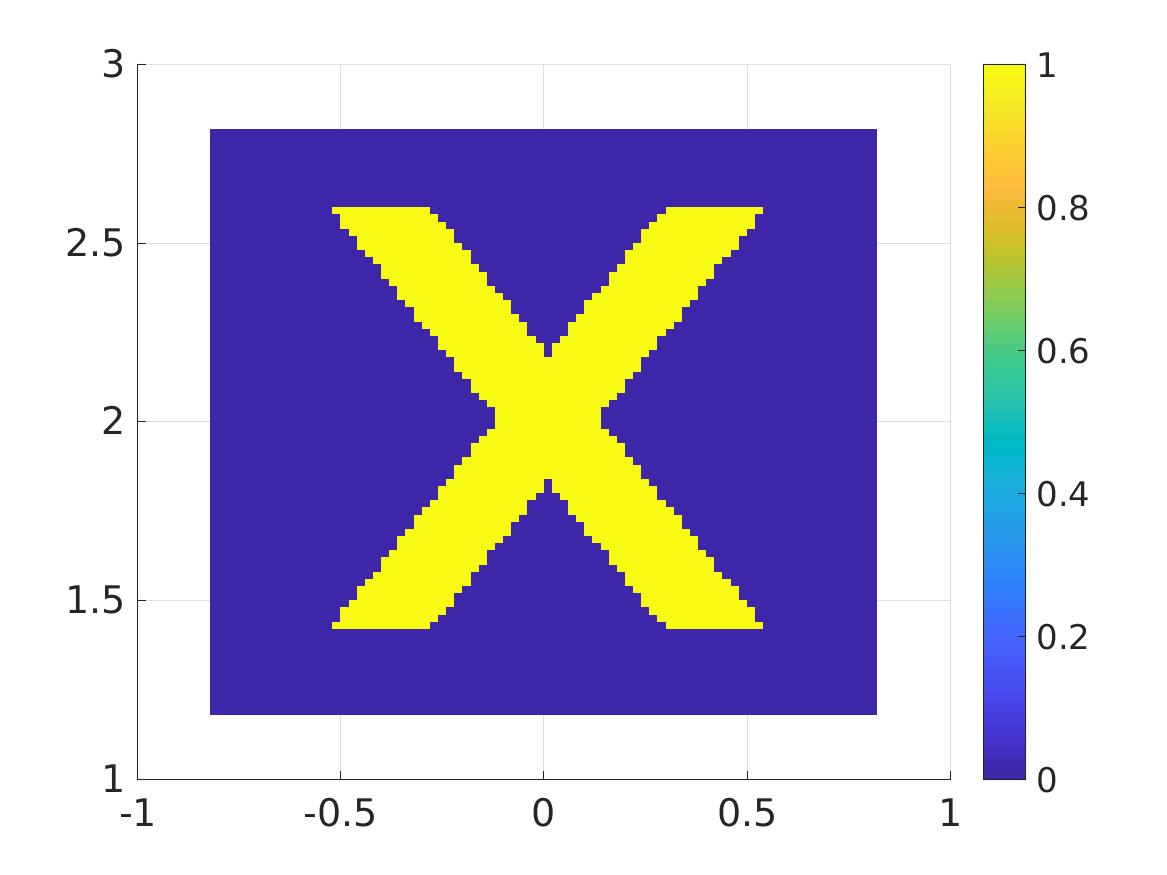}}\quad 
\subfloat[The reconstructed 
function $f^{\mathrm{comp}}$, no noise]{\includegraphics[width
=.3\textwidth]{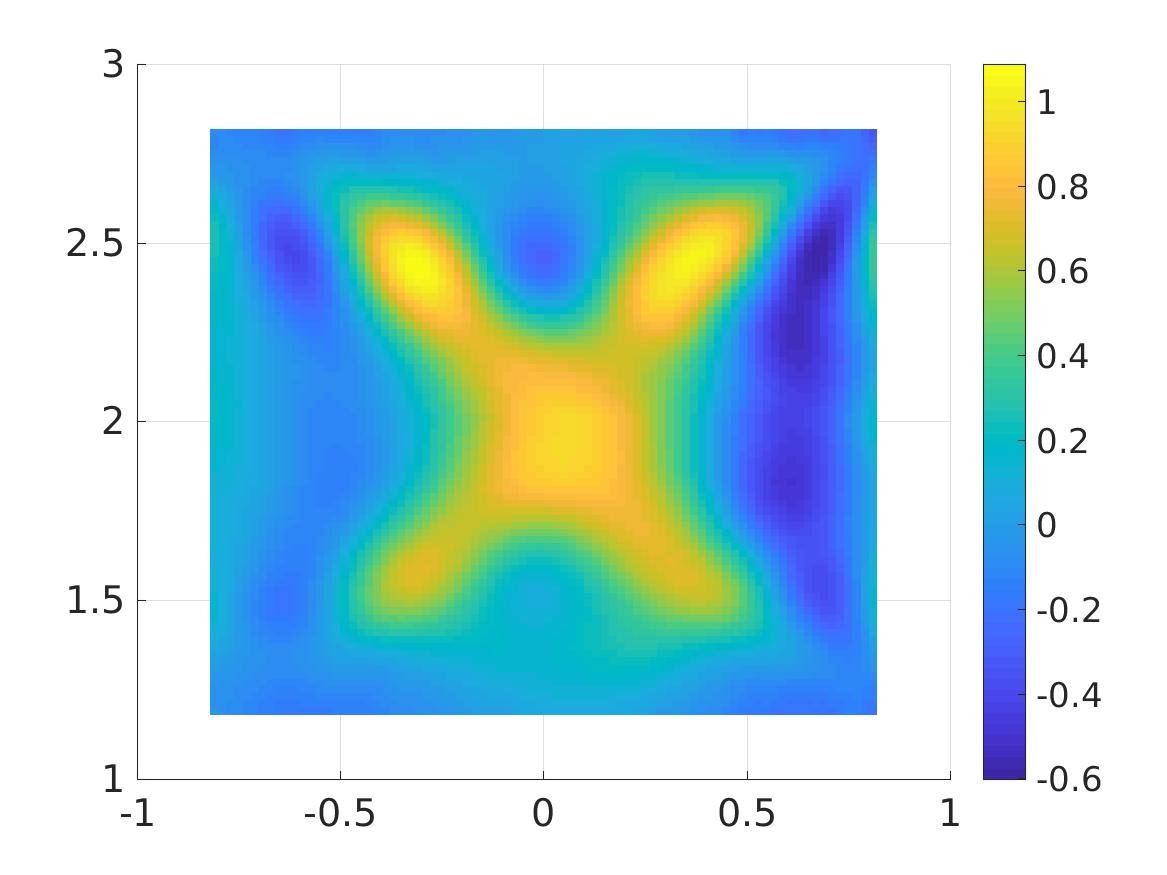}} 
\subfloat[The reconstructed 
function $f^{\mathrm{comp}}_{\delta}$, noise 30\%]{\includegraphics[width =
.3\textwidth]{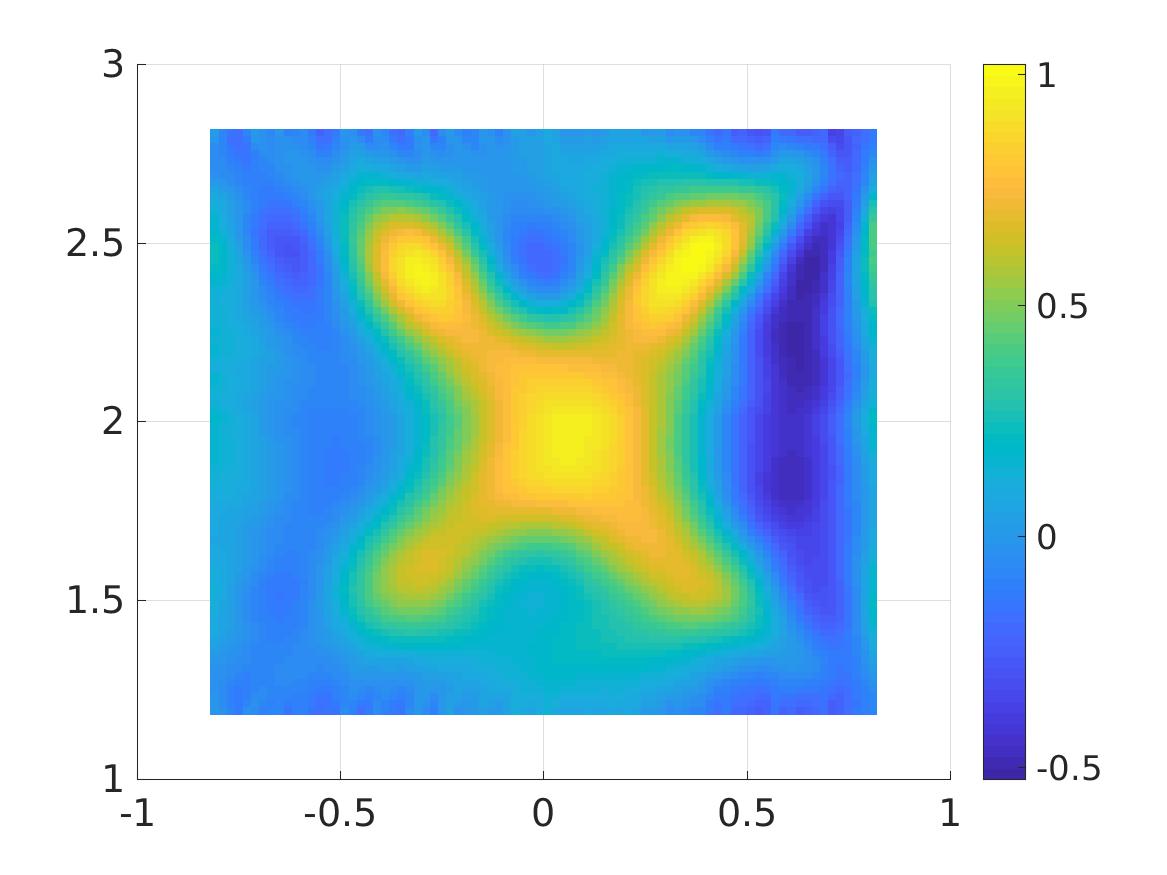}}
\par
\subfloat[\label{incomplete2}The incomplete tomographic
data, no noise]{\includegraphics[width =
.3\textwidth]{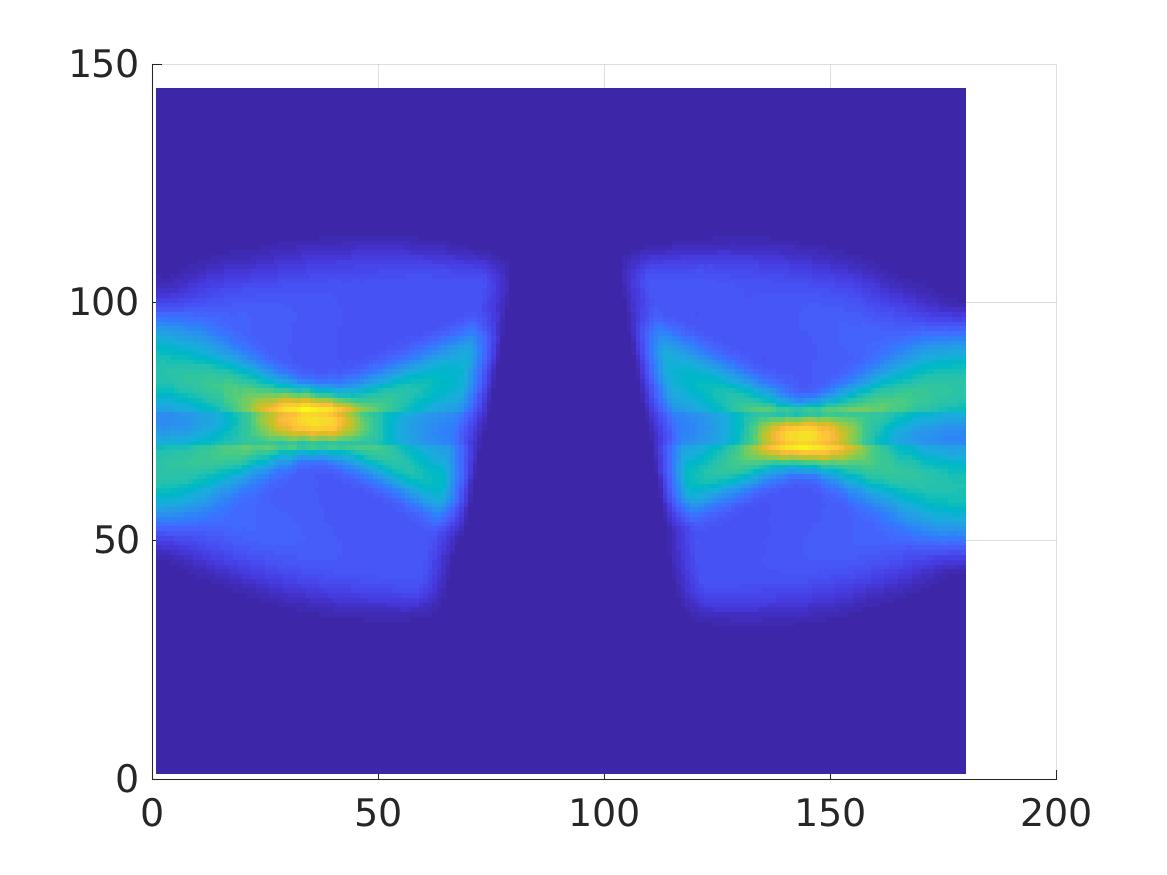}} \quad 
\subfloat[Post-processed $f^{\mathrm{comp}}$,  no
noise]{\includegraphics[width = .3\textwidth]{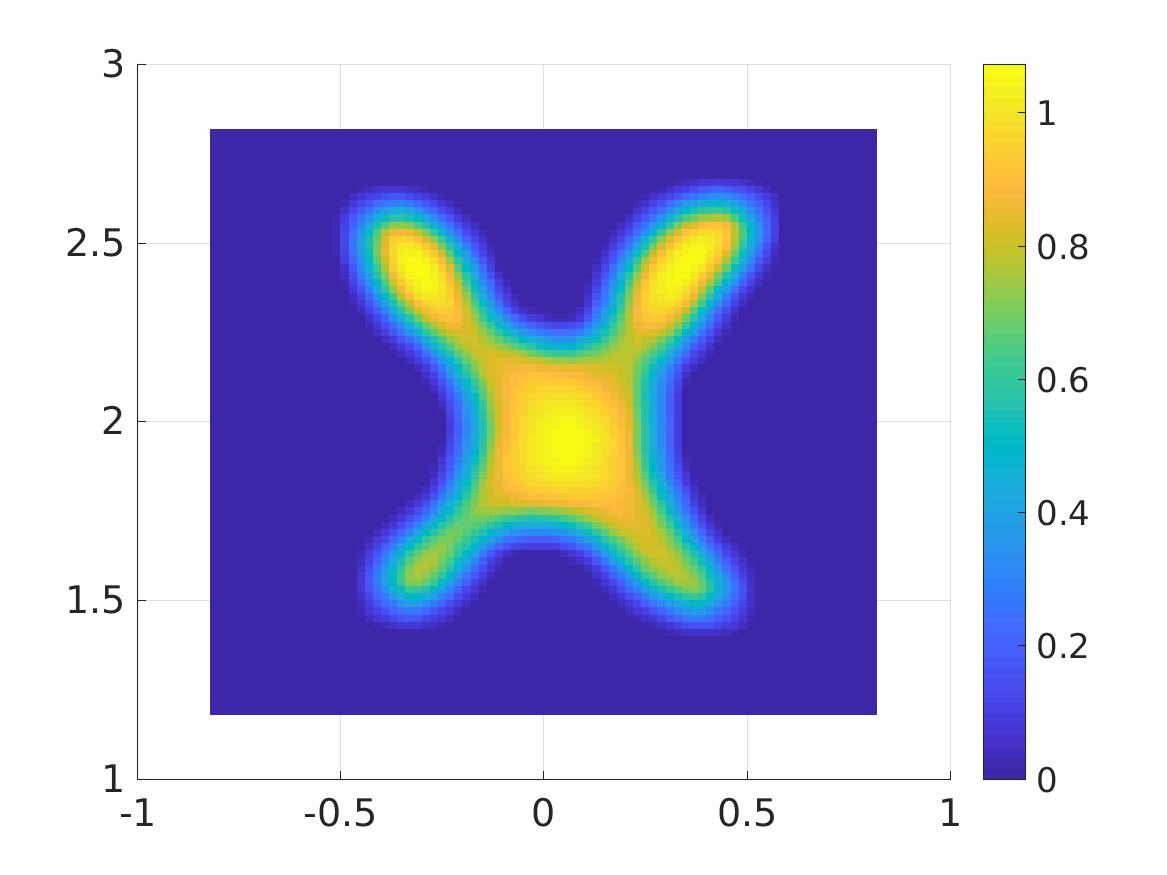}} 
\subfloat[Post-processed $f^{\mathrm{comp}}_{\delta}$, noise level 30\%]{\includegraphics[width =
.3\textwidth]{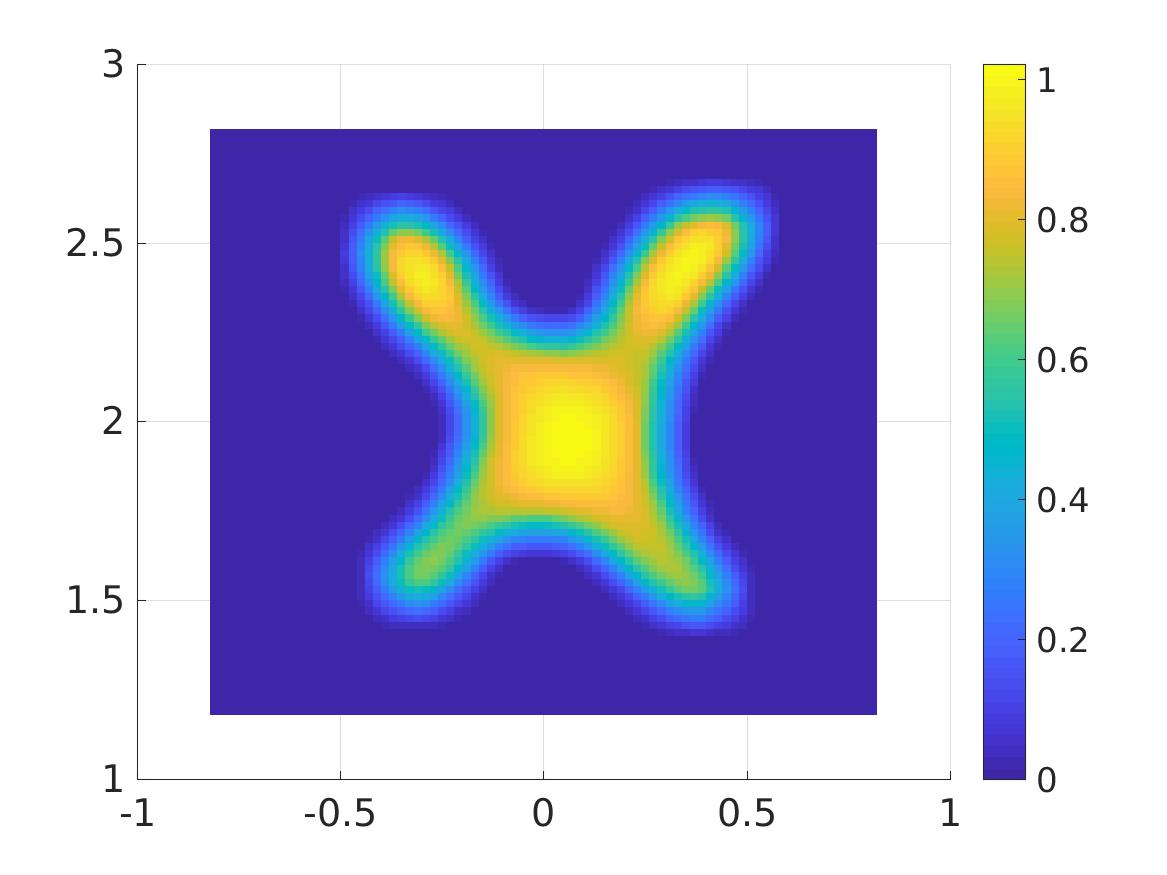}}
\end{center}
\caption{\textit{Test 2. The true and reconstructed source functions for for X-shaped non-smooth inclusion, for the case of uniform scattering}}
\label{fig model 2}
\end{figure}

\item \textbf{Test 3. Y-shaped inclusion, strongly forward-peaked scattering}

The function $f_{\mathrm{true}}$ depicted on Figure \ref{fig model 3}(a) is
smooth and compactly supported in the domain. We define absorption and
scattering coefficients as 
\begin{equation*}
\mu _{a}(\mathbf{x})=\left\{ 
\begin{array}{rl}
0.15 & \mbox{if }(x,y)\in \mathrm{supp}(f_{\mathrm{true}}), \\ 
0.1 & \mbox{if }(x^{2}+y^{2})<0.8, \\ 
0 & \mbox{otherwise,}%
\end{array}%
\right.
\end{equation*}
and
\begin{equation*}
 \mu _{s}(\mathbf{x})=\left\{ 
\begin{array}{rl}
0.01 & \mbox{if }(x^{2}+y^{2})<0.8, \\ 
0 & \mbox{otherwise,}%
\end{array}%
\right.
\end{equation*}%
and the scattering phase function 
\begin{equation}
K(\mathbf{x},\alpha ,\beta )=H(\alpha ,\beta )=\frac{1}{2d}\left[ \frac{%
1-g^{2}}{1+g^{2}-2g\cos {(\alpha -\beta })}\right],
\label{6.6}
\end{equation}%
where $H(\alpha ,\beta )$ is the 2-dimensional Henyey-Greenstein function. In (\ref{6.6}), $g=g(x,y)$ is the Henyey-Greenstein factor, it is
the smoothed version of the function 
\begin{equation*}
\tilde{g}(x,y)=\left\{ 
\begin{array}{rl}
0.9 & \mbox{if }(x^{2}+y^{2})<0.8, \\ 
0.5 & \mbox{otherwise,}%
\end{array}%
\right. .
\end{equation*}%
The numerical solution for Test 3 is depicted on Figure \ref{fig model 3}.
\end{enumerate}

\begin{figure}
\begin{center}
\subfloat[The function $f_{\rm true}$]{\includegraphics[width =
.3\textwidth]{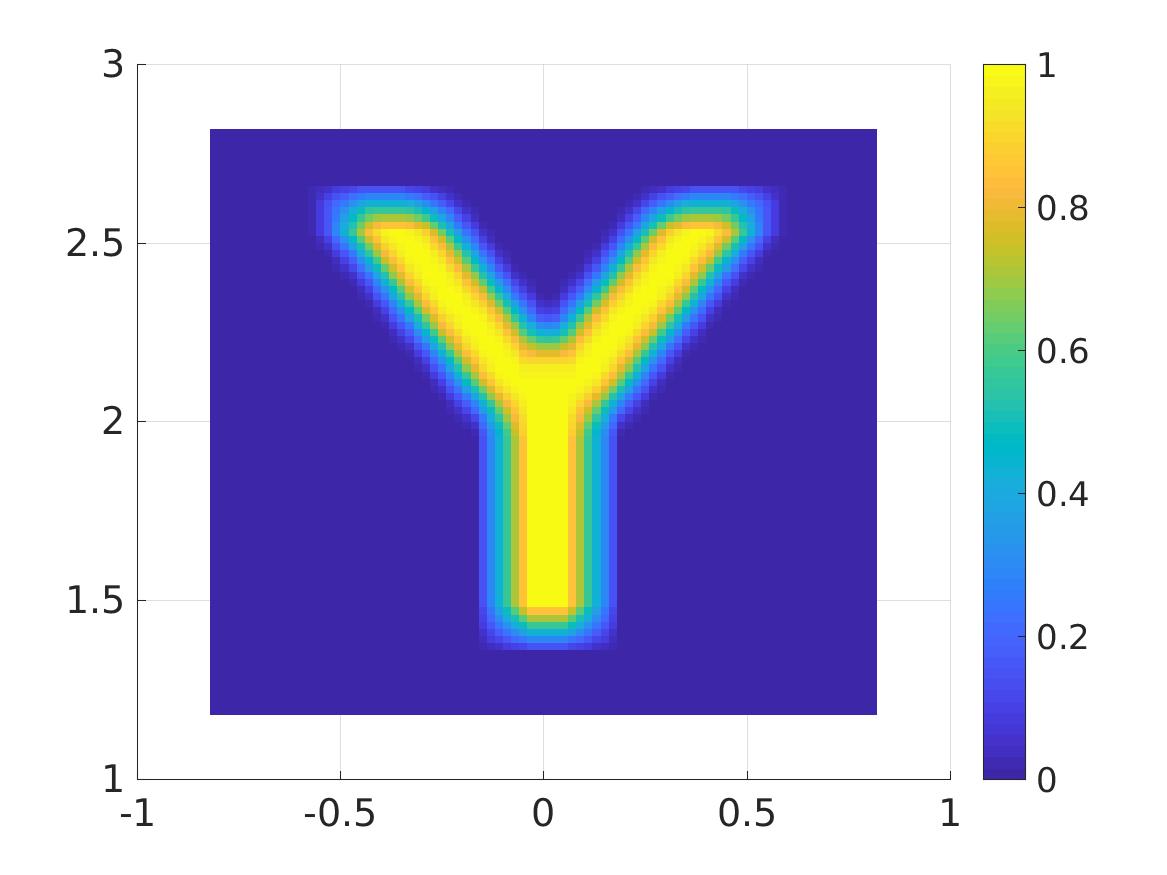}}\quad 
\subfloat[The reconstructed 
function $f^{\mathrm{comp}}_{\delta}$, noise 60\%]{\includegraphics[width =
.3\textwidth]{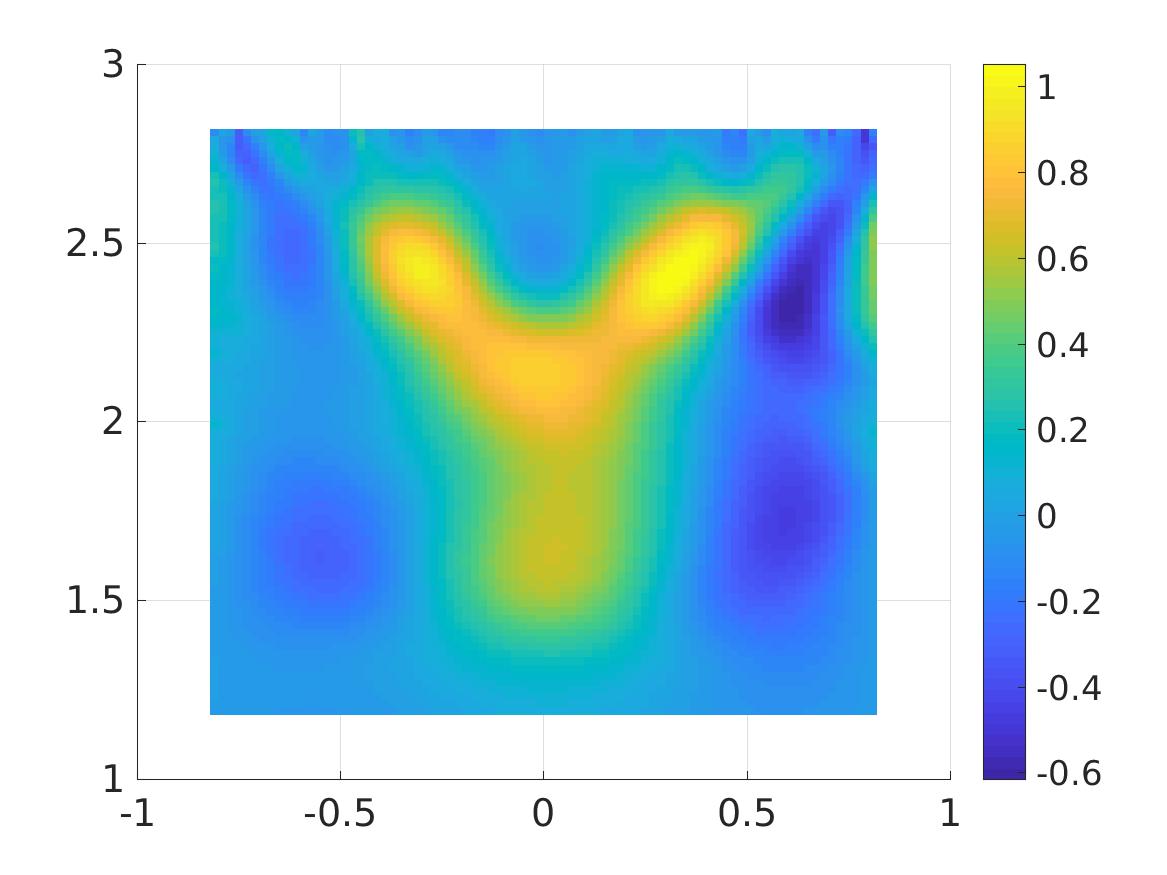}} 
\subfloat[Post-processed $f^{\mathrm{comp}}_{\delta}$,  noise level 60\%]{\includegraphics[width =
.3\textwidth]{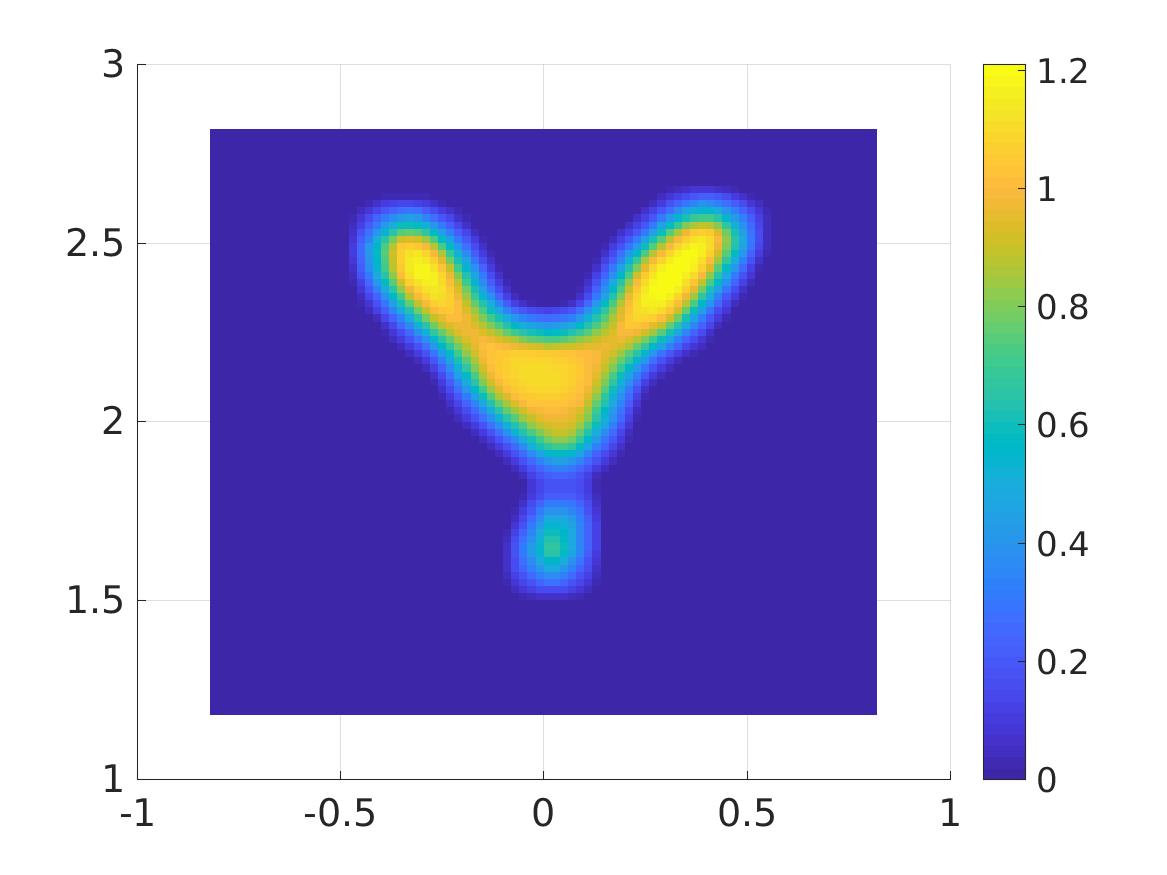}}
\par
\subfloat[Post-processed $f^{\mathrm{comp}}$,  no noise]{\includegraphics[width =
.3\textwidth]{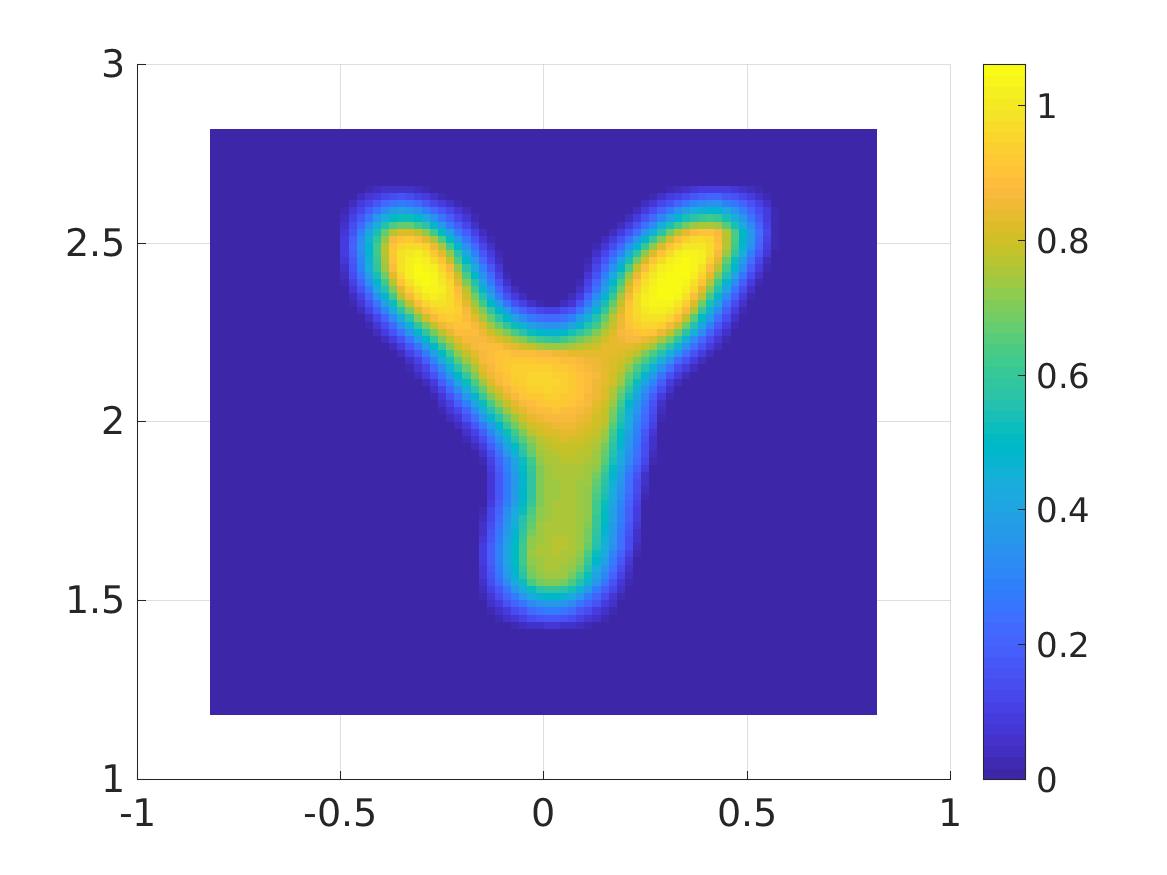}} 
\subfloat[The absorption coefficient]{\includegraphics[width =
.3\textwidth]{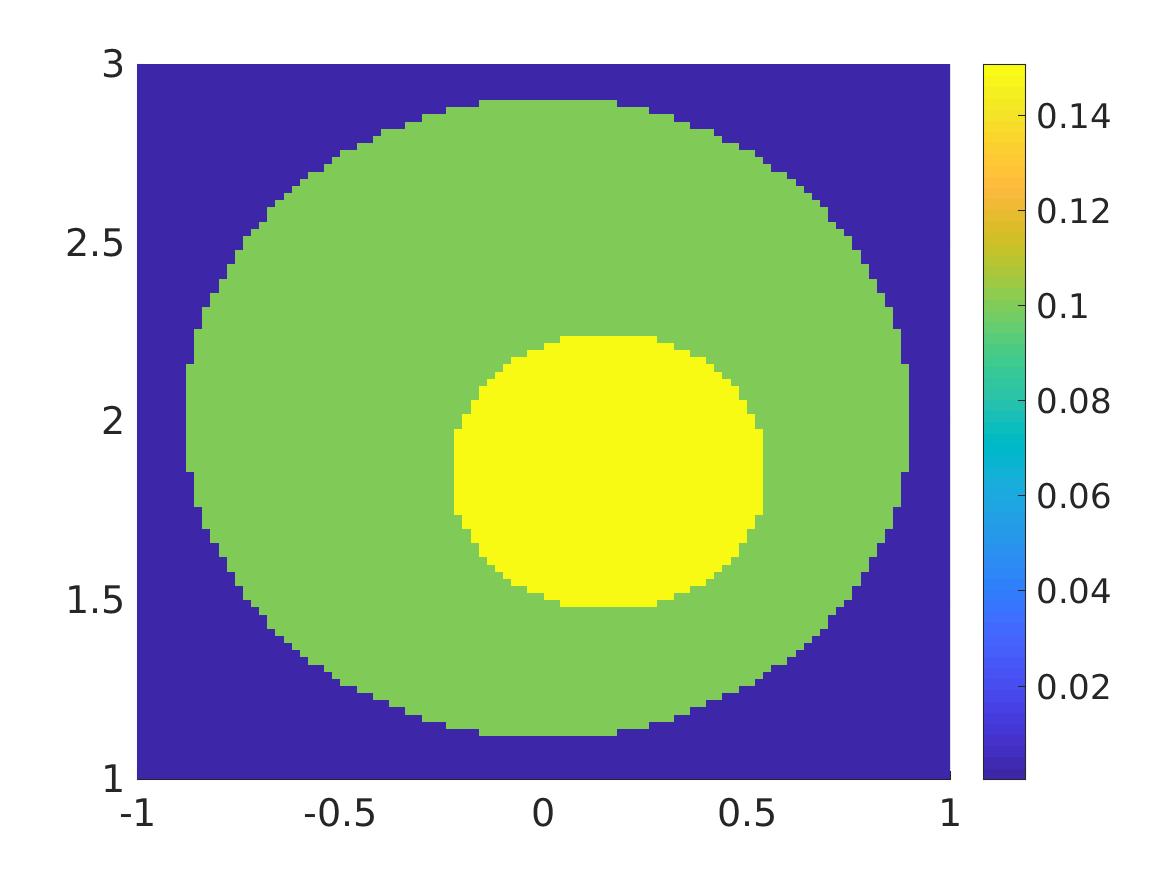}} 
\subfloat[The scattering coefficient]{\includegraphics[width =
.3\textwidth]{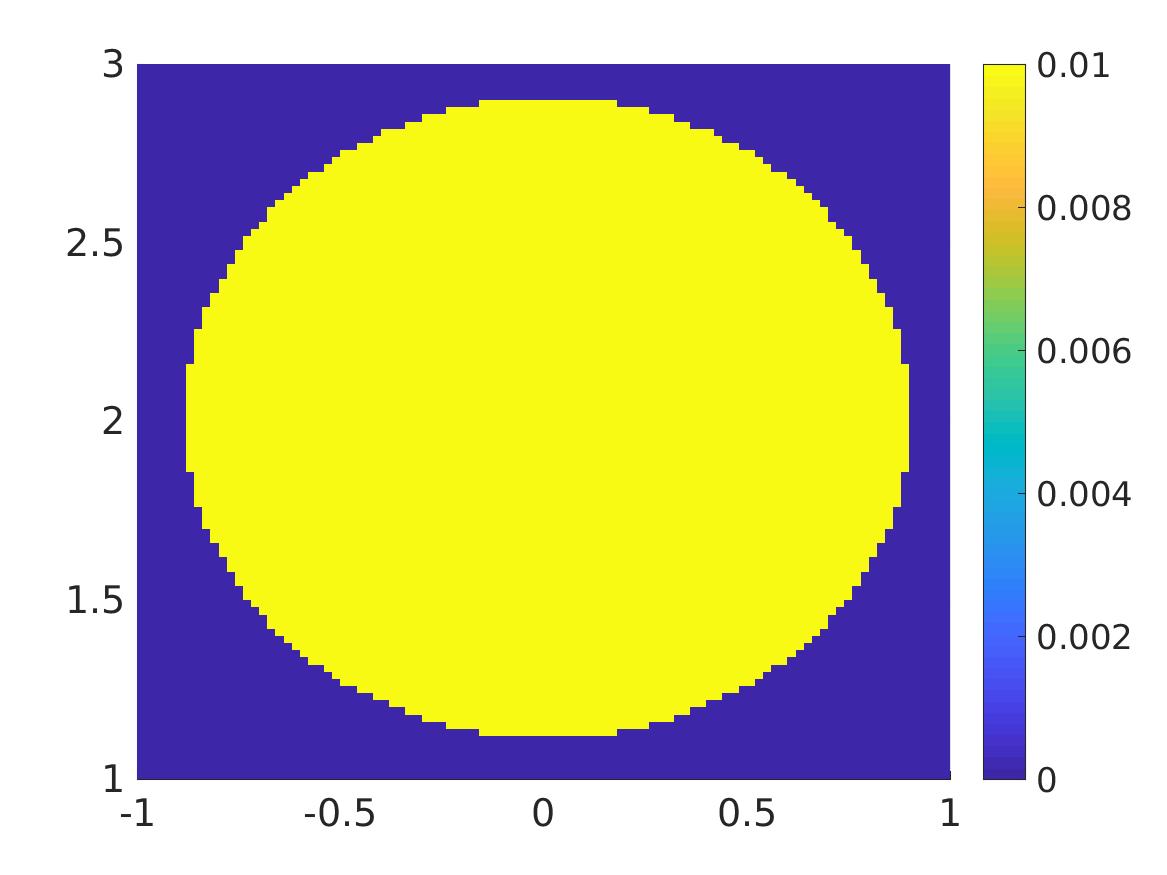}}
\end{center}
\caption{\textit{Test 3. The true and reconstructed source functions for
Y-shaped inclusion for the case of strongly forward-peaked scattering}}
\label{fig model 3}
\end{figure}

\begin{remark}
 In Tests 1 and 2, Figures \ref{incomplete1} and \ref{incomplete2} are given only
for the illustration purpose. They represent the incomplete Radon transform
data, obtained by the well-known Radon transform of the function $f^{true}$,
in which case $\mu _{a}=\mu _{s}=0.$ The case when the data is available for
all $\mathbf{x}$ and $\mathbf{x}_{\alpha }$ such that the set of lines $L(%
\mathbf{x},\mathbf{x}_{\alpha })$ contains all possible lines intersecting
domain $\Omega $ is considered to be the tomographic inverse problem with
complete data. In contrast to this, we study the case when, the point source 
$\mathbf{x}_{\alpha }$ runs along the straight line as shown in Figure \ref%
{fig0}. In that scenario, the data in our problem is said to be incomplete
and Figures \ref{incomplete1} and \ref{incomplete2} illustrate this. In Test 3 we omit the image
representing incomplete Radon data, for the reason that it does not differ
qualitatively from  Figures \ref{incomplete1} and \ref{incomplete2}, while the distributions of the
absorption and scattering coefficients differ significantly from the
previous two tests.
\end{remark}

\section{Appendix}

\begin{theorem}[Uniqueness and existence of the solution of the Forward Problem]
Consider a rectangular domain $\Omega $, defined in (\ref{2.1}). Assume that
functions $f, \mu_a, \mu_s \in C^{1}\left( \mathbb{R}^{2}\right), f\in L^{2}\left( \mathbb{R}^{2}\right) $ and the function $K\in C^{1}\left( 
\mathbb{R}^{2}\times \left[ -d,d\right] ^{2}\right) $. Also, assume that (%
\ref{2.4}) holds. Then there exists a unique solution $u(\mathbf{x},\alpha )$
of Problem \ref{pro forward} in the domain $\Omega $ such that $u\left( x,y,\alpha
\right) \in C^{1}\left( \mathbb{R}\times \left[ 0,b\right] \times \left[ -d,d%
\right] \right) .$ Furthermore, $u(x,y,\alpha )=0$ for $y\in (0,a)$ as well
as for sufficiently large $\left\vert x\right\vert $. \label{thm3}
\end{theorem}

\textbf{Proof.} Let $L(\mathbf{x},\mathbf{x}_{\alpha })$ be the segment of
the straight line connecting points $\mathbf{x}$ and $\mathbf{x}_{\alpha }$.
Then for any appropriate function $\phi (\mathbf{x})$ 
\begin{equation*}
\int \displaylimits_{L(\mathbf{x},\mathbf{x}_{\alpha })}\phi (\boldsymbol{%
\xi })ds_{\boldsymbol{\xi }}=\frac{|\mathbf{x}-\mathbf{x}_{\alpha }|}{y}%
\int_{0}^{y}\phi \left( \alpha +\frac{w(x-\alpha )}{y},w\right) dw,
\end{equation*}%
where $ds_{\boldsymbol{\xi }}$ is the arc length. It follows from (\ref{2.1}%
)-(\ref{2.4}) that $u(\mathbf{x},\alpha )=0$ for $y\in (0,a).$ Therefore, the
Forward Problem in the domain $\Omega $ is equivalent to 
\begin{multline}
u(\mathbf{x},\alpha ) =\chi ^{-1}(\mathbf{x},\alpha )\frac{|\mathbf{x}-%
\mathbf{x}_{\alpha }|}{y}\int \displaylimits_{a}^{y}\chi (\mathbf{z})\mu
_{s}(\mathbf{z})\int \displaylimits_{-d}^{d}K(\mathbf{z},\alpha ,\beta )u(%
\mathbf{z},\beta )d\beta dw
\\ +\chi ^{-1}(\mathbf{x},\alpha )\frac{|\mathbf{x}-\mathbf{x}_{\alpha }|}{y}%
\int \displaylimits_{a}^{y}(f\chi )(\mathbf{z})dw,\hspace{0.5em}
\label{7.2}
\end{multline}
and
\begin{equation}
\mathbf{z}(w,x,\alpha )=\left( \alpha +\frac{w(x-\alpha )}{y},w\right)
,\quad \ln {\chi (\mathbf{x},\alpha )}={\frac{|\mathbf{x}-\mathbf{x}_{\alpha
}|}{y}\int_{0}^{y}}(\mu _{a}+\mu _{s})(\mathbf{z})dw.  \label{7.3}
\end{equation}%
Estimate from the below the absolute value of the first argument in $\mathbf{%
z}(w,x,\alpha )$ in (\ref{7.3}). By (\ref{2.1}), (\ref{2.4}) and and (\ref%
{7.3}) the left hand side of equation (\ref{7.2}) is not zero only if $y\in
\left( a,b\right) .$ Since in $\mathbf{z}(w,x,\alpha )$ we have $w\in \left(
a,y\right) ,\alpha \in \left( -d,d\right) ,$ then%
\begin{equation}
\left\vert \alpha +\frac{w(x-\alpha )}{y}\right\vert \geq \frac{w}{y}%
\left\vert x-\alpha \right\vert -\left\vert \alpha \right\vert \geq \frac{a}{%
b}\left\vert x\right\vert -\left( 1+\frac{a}{b}\right) \left\vert \alpha
\right\vert \geq \frac{a}{b}\left\vert x\right\vert -\left( 1+\frac{a}{b}%
\right) d.  \label{7.4}
\end{equation}%
Suppose that $\left\vert x\right\vert \geq X$ and $X$ is so large that 
$
1-\frac{1}{X}\left( 1+\frac{b}{a}\right) d>\frac{1}{2},\text{ }\frac{a}{2b}%
X>R.
$
Then (\ref{7.4}) implies that $\left\vert \alpha +w(x-\alpha )/y\right\vert
>R.$ Hence, by (\ref{2.1}) and (\ref{7.3}) the right hand of equation (\ref%
{7.2}) equals zero for $\left\vert x\right\vert \geq X.$ Let 
\begin{equation*}
A\left( u\right) :C\left( \mathbb{R}^{2}\times \left[ -d,d\right] \right)
\rightarrow C\left( \left( \left\vert x\right\vert \leq X\right) \times %
\left[ -d,d\right] \times \left[ -d,d\right] \right)
\end{equation*}%
be the operator in the right hand side of (\ref{7.2}). Then (\ref{7.2}) can
be considered as the equation $u=A\left( u\right) $ with the Volterra-like
integral operator, where the \textquotedblleft Volterra property" is due to
the integration with respect to $y$. Therefore, the latter equation can be
solved iteratively, as it is usually done for the Volterra integral
equations. It is obvious from the above discussion that all iterates $%
A\left( u_{n}\right) \left( \mathbf{x},\alpha \right) =0$ for $\left\vert
x\right\vert \geq X.$ Thus, the solution of equation (\ref{7.2}) in the
space $C\left( \mathbb{R}^{2}\times \left[ -d,d\right] \right) $ exists and
is unique. The $C^{1}-$smoothness of this solution with respect to $%
x,y,\alpha $ obviously follows from the well known convergence estimate for
the iterates of a Volterra integral equation. Due to the above mentioned
equivalence, this implies uniqueness and existence of the solution of the
Forward Problem $\ \square $

The numerical solution of the Forward Problem was performed via the
iterative solution of the Volterra-like integral equation (\ref{7.2}) for $(%
\mathbf{x},\alpha )\in \left( \left\vert x\right\vert \leq X\right) \times %
\left[ a,b\right] \times \left[ -d,d\right] $.


\end{document}